\documentclass[11pt,epsfig]{article}
\usepackage{mathrsfs}
\usepackage{amsmath}
\usepackage{amsthm}
\usepackage{epsfig}
\usepackage{hyperref}
\usepackage{amsfonts}
\usepackage{color}
\usepackage[margin=1.1in]{geometry}

\newtheoremstyle{thmm}{1.5ex plus 1ex minus .2ex}{1.5ex plus 1ex minus
.2ex}{\rmfamily}{}{\bfseries}{}{1em}{} \theoremstyle{thmm}
\newtheorem{theorem}{Theorem}[section]
\newtheorem{lemma}{Lemma}[section]

\renewenvironment{proof}[1][Proof]{\noindent\textit{#1. } }{\hfill$\square$=
}

\allowdisplaybreaks

\renewcommand{\theequation}{\thesection.\arabic{equation}}
\newcommand{\nn}{\nonumber}
\def \endproof{\vrule height8pt width 5pt depth
0pt}
\def\refe#1{(\ref{#1})}

\def\d{\delta}

\def\R{\mathbb{R}}

\def\d{{\rm d}}

\begin{document}
\title{\bf Linearized
FE approximations to a strongly nonlinear diffusion equation 
\\[10pt]} 
\author{Buyang
Li\,\setcounter{footnote}{0}\footnote{Department
of Mathematics, Nanjing University, Nanjing 210093, Jiangsu, P.R. China.
The work of the author was supported in part by a grant from NSFC
(Grant No. 11301262)
{\tt buyangli@nju.edu.cn}. }
~~and ~Weiwei Sun\footnote{Department
of Mathematics, City University of Hong Kong,
Kowloon, Hong Kong.  The work of the
author was supported in part by a grant from the Research Grants
Council of the Hong Kong Special Administrative Region, China
(Project No. CityU 102005)  {\tt maweiw@math.cityu.edu.hk}. }}

\date{}
\maketitle

\begin{abstract}
We study fully discrete linearized Galerkin finite element 
approximations to a nonlinear gradient flow, 
applications of which can be found in many  
areas. Due to the strong nonlinearity of the equation, 
existing analyses for implicit schemes require certain restrictions on the time step 
and no analysis has been explored for linearized schemes.
This paper focuses on the unconditionally optimal 
$L^2$ error estimate of  a linearized scheme. 
The key to our analysis is 
an iterated sequence of time-discrete elliptic equations and
a rigorous analysis of its solution.
We prove the $W^{1,\infty}$ boundedness of 
the solution of the time-discrete system
and the corresponding finite element solution, based on a more precise estimate 
of elliptic PDEs 
in $W^{2,2+\epsilon_1}$ and $H^{2+\epsilon_2}$ 
and a physical feature of the gradient-dependent diffusion 
coefficient. Numerical examples are provided to support our theoretical analysis.
\\

\medskip

\end{abstract}

~~~{{\bf AMS subject classifications:} 65N12,
65N30, 35K61.}

\section{Introduction}
\setcounter{equation}{0}
We consider the nonlinear diffusion equation
\begin{align}
&\frac{\partial u}{\partial t}-\nabla\cdot(\sigma(|\nabla
u|^2)\nabla u)=g \label{e-parab-1}
\end{align}
in a convex polygonal domain $\Omega$ in $\R^2$
with the Neumman boundary condition
\begin{align}
\label{BC}
\begin{array}{ll}
\nabla u\cdot{\vec n}=0~~ &\mbox{on}~~\partial\Omega
\end{array}
\end{align}
and the initial condition
\begin{align}
\label{IniC}
\begin{array}{ll}
u(x,0)=u_0(x)~~
&\mbox{for}~~x\in\Omega ,
\end{array}
\end{align}
where $g$ is a given function and 
\begin{align} \label{sigma}
\sigma(s^2) = \frac{1}{\sqrt{\lambda^2 + s^2}}  
\end{align}
is a gradient-dependent
diffusion coefficient, where $\lambda$ is a positive constant. The equation has been involved in many applications,
such as minimal surface flow \cite{LT}, prescribed mean curvature flow \cite{DD,Ger},
geometric measure theory \cite{BC}, and a regularized model in 
image denoising \cite{CL,CX,CC,EV,GMS,LC,LG,PM,VO}.
A review article for the applications in image processing 
was given in \cite{CMY}.  

Mathematical analysis of the nonlinear diffusion equation 
\refe{e-parab-1} was
studied in \cite{FP,Ger}. In particular, the $W^{1,\infty}$ regularity of the 
solution was
proved in \cite{FP}, which implies arbitrarily 
higher regularity of the solution in a smooth domain 
(by the method of Section 8.3.2 of \cite{Evans}).
Numerical methods and simulations for the nonlinear
diffusion equation have been investigated extensively
in the last several dacades. For examples, see
\cite{ABS,LG,LC,VO} for finite difference methods and 
 \cite{CX,DD,DER,EV,FOP,FP,FP2} for finite element methods (FEMs).
Explicit schemes may not be efficient due to their strong time-step restrictions.
A fully implicit backward Euler--Galerkin FEM was analyzed in \cite{FP},
where optimal convergence rate was proved under the condition $\tau = O(h^2)$.
Suboptimal error estimates for the scheme were presented in
\cite{FP2} under a weaker mesh restriction $\tau=o(h^{1/2})$, and
further analysis on the convergence rate of the scheme
with respect to the regularization parameter was given in \cite{FOP}.
The implicit backward Euler scheme was also studied in \cite{EV}
with a lumped mass FEM, where $L^\infty$-boundedness
of the numerical solution was proved and no error estimates were presented.
In these fully implicit schemes,
one has to solve a system of nonlinear equations at each time step and
an extra inner iteration is needed. 
In addition to the implicit schemes, linearized semi-implicit
FEMs for the nonlinear diffusion equation have also been investigated by
several authors \cite{CX,LC,MS}.  In this method, the gradient-dependent
diffusion coefficient is calculated with the numerical solution at
the last time step and Galerkin FEMs are used to solve the
linearized equation. The scheme only requires the solution of
a linear system at each time step, which is simple and efficient
for implementation.  However, theoretical error analysis 
of the linearized scheme
seems very difficult due to the strong nonlinear structure.
As far as we know, no optimal error estimates of linearized semi-implicit FEMs
are available for the nonlinear diffusion equation.
The major difficulty for the analysis of the semi-implicit scheme is due to
the nature of the linearization of the scheme, which leads to the arising of
the energy-norm errors at two different time levels in the error equation (see \refe{errFEMFDFEM}-\refe{dkhf2} for the estimates of the error equation).

In this paper,  we study linearized backward Euler--Galerkin methods 
for the nonlinear system
\refe{e-parab-1}-\refe{IniC}. Our focus is on unconditionally optimal error 
estimates of numerical methods. The key issue in the analysis is  
to establish the $W^{1,\infty}$ convergence of the numerical solution.
To deal with the strong nonlinearity from the gradient-dependent
diffusion coefficient, we introduce an iterated sequence of time-discrete 
elliptic PDEs as in \cite{LS1,LS2}. 
Thus the linearized backward Euler--Galerkin method coincides with
the corresponding FE approximation to the time-discrete system.
We prove the $W^{1,\infty}$ convergence of the solution of the time-discrete system 
and FE solution, in terms of 
a more precise estimate for elliptic PDEs in $W^{2,2+\epsilon_1}$ and 
$H^{2+\epsilon_2}$:
\begin{align*} 
& \| u \|_{L^{2+\epsilon_1}} \le (1+ \epsilon^*_1) \| \Delta u \|_{L^{2+\epsilon_1}} 
\\ 
& \| u \|_{H^{2+\epsilon_2}} \le (1+ \epsilon^*_2) \| \Delta u \|_{H^{\epsilon_2}}  
\, , 
\end{align*} 
and a physical feature of the gradient-dependent diffusion coefficient: 
$2|\sigma'(s^2)|s^2<\sigma(s^2)$. With these a priori estimates, we establish 
the $L^2$-norm optimal error estimate without any time-step restrictions. 

The rest part of this paper is organized as follows.  In Section 2, we
introduce some notations and the linearized
backward Euler--Galerkin FEM for the nonlinear diffusion
equation \refe{e-parab-1}-\refe{IniC}, and then we present our
main results and our methodology.
In Section 3, we prove our main results based on the regularity
and $W^{1,\infty}$-convergence of the time-discrete solution, 
while the rigorous proof of the regularity and $W^{1,\infty}$-convergence 
of the time-discrete solution is postponed to Section 4.
Numerical examples are presented in Section 5, which confirm our
theoretical analysis and show clearly that the linearized scheme is efficient and 
no time-step conditions are needed.

\section{Notations and main results}
\setcounter{equation}{0}
Let $\Omega$ be a given convex polygon in $\R^2$.
For $1\leq p\leq \infty$ and
any nonnegative integer $k$, we denote by $W^{k,p}(\Omega)$ the usual Sobolev 
space of functions defined on $\Omega$ and, to simplify the notations, 
we set $W^{k,p}:=W^{k,p}(\Omega)$, $H^k:=W^{k,2}(\Omega)$ and $L^p:=W^{0,p}$. 
For $s\in(0,1)$, we define $H^{k+s}:=(H^{k},H^{k+1})_{[s]}$ as the complex 
interpolation space between $H^k$ and $H^{k+1}$. 
More detailed discussions for the complex interpolation spaces can be found in 
literature, $e.g.$, see the classical book \cite{BL} by Bergh and L\"ofstr\"om. 

For a given quasi-uniform triangulation of
$\Omega$ into triangles $T_j$,
$j=1,\cdots,J$, we denote by $h=\max_{1\leq j\leq J}\{\mbox{diam}\,T_j\}$
the mesh size and define a finite element space by
\begin{align*}
&V_h^r=\{v_h\in C(\overline\Omega): v_h|_{T_j}
\mbox{~is~a~polynomial of degree $r$}\}
\end{align*}
so that $V_h^r$ is a subspace of $H^1(\Omega)$.
Let $\Pi_h:C(\overline\Omega)\rightarrow V_h^r$ denote the Lagrangian
interpolation operator. Let
$0=t_0<t_1<\cdots<t_N=T$ be a uniform partition of the time interval
$[0,T]$ with $t_n=n\tau$.  For a sequence of functions $\{ f^n \}_{n=0}^N$,
we define a time-difference operator by
\begin{align}
& D_\tau f^{n+1}=\frac{f^{n+1}-f^n}{\tau}, \quad\mbox{for
$n=0,1,\cdots,N-1$.} \label{dfd}
\end{align}


We define the linearized backward Euler--Galerkin finite element scheme by
\begin{align}
&\big(D_\tau U^{n+1}_h,v\big)+\big(\sigma(|\nabla U^n_h|^2)
\nabla
U^{n+1}_h,\nabla v\big) =\big(g^{n+1},v\big),\quad \forall
v\in V_h^r
,\label{FDFEM}
\end{align}
with the initial condition $U^0_h= \Pi_hu_0$ and $r\ge 2$.
At each time step, the scheme only requires the solution of a linear system.
Also we assume   that the solution of {\rm
(\ref{e-parab-1})-(\ref{IniC})}
exists and satisfies 
\begin{align}
\label{regularity} 
&\|u_0\|_{H^{r+1}}+\|u\|_{L^\infty((0,T);H^{r+1})}
+\|\partial_{t}u\|_{L^\infty((0,T);H^{r+1})}
+\|\partial_{tt}u\|_{L^2((0,T);L^2)}  \leq M_0,
\end{align}
where $M_0$ is some positive constant. For simplicity, 
we assume that $g=g(x,t)$ in this paper. The analysis presented in this paper can be easily extended to the general case $g=g(u,x,t)$ for the scheme
\begin{align*}
&\big(D_\tau U^{n+1}_h,v\big)+\big(\sigma(|\nabla U^n_h|^2)
\nabla
U^{n+1}_h,\nabla v\big) =\big(g(U^n_h,x,t^{n}),v\big),\quad \forall
v\in V_h^r
,
\end{align*}
if $g$ is a smooth function of $u$, $x$ and $t$.

Our main results are given in the following theorem concerning
the unconditionally optimal convergence rate of
the numerical solution.
\begin{theorem}
\label{ErrestFEMSol}
{\it
Suppose that the system {\rm
\refe{e-parab-1}-\refe{IniC}} has a unique solution $u$
satisfying
the regularity condition {\rm \refe{regularity}}.
Then there exists a positive constant $C_0$, independent of
$\tau$ and
$h$, such that the finite element system {\rm
(\ref{FDFEM})}
admits a unique solution $\{ U^n_h \}_{n=1}^N$ satisfying
\begin{align}\label{optimalL2est}
&\|U_h^n-u^n\|_{L^2}\leq C_0(\tau+h^{r+1}) .
\end{align}
}
\end{theorem}
\medskip

To prove the above theorem,
we introduce an iterated sequence of elliptic PDEs (time-discrete system) as
proposed in \cite{LS1,LS2}:
\begin{align}
&D_\tau U^{n+1}-\nabla\cdot(\sigma(|\nabla U^n|^2)\nabla
U^{n+1})=g^{n+1} \label{e-TDparab-1},
\end{align}
with the boundary condition $\nabla U^{n+1}\cdot{\vec n}=0$ on
$\partial\Omega$ and the initial condition $U^0=u_0$.
Then the fully
discrete solution $U^{n+1}_h$ coincides with the finite element solution
of (\ref{e-TDparab-1}). In view of this property, we split the error into
$$
U^n_h-u^n =(U^n_h-U^n)+(U^n-u^n)
$$
and analyze the two error functions separately.
The regularity of the solution of the time-discrete
system (\ref{e-TDparab-1}) is given in the following theorem.

\begin{theorem}\label{ErrestDisSol}
{\it
Under the assumption of Theorem {\rm\ref{ErrestFEMSol}},
there exist positive constants $\tau_0^*$, $C_0^*$, $p>2$ and $s_0>0$, which
are dependent only on  $M_0$, $\Omega$ and $T$ and independent
of $\tau$ and $h$, such that when $\tau<\tau_0^*$
the time-discrete system
{\rm(\ref{e-TDparab-1})} admits a unique solution $\{ U^{n}
\}_{n=0}^N$ satisfying
\begin{align}
& \max_{0 \le n \le N} (\|U^n\|_{W^{2,p}}^2 +\|U^n\|_{H^{2+s_0}}^2)+
\sum_{n=1}^N\tau\|
D_\tau
U^{n}\|_{H^2}^2  \leq C_0^*, \label{dt-bound}\\
&\max_{1\leq n \leq N}
\|e^n\|_{H^1}^2+\sum_{n=1}^N\tau\| e^{n}\|_{H^2}^2
+\sum_{n=1}^N\tau\| D_\tau e^{n}\|_{L^2}^2\leq
C_0^*\tau^2, \label{dt-b33}\\
& \max_{1\leq n \leq N}\|e^n\|_{W^{2,p}}\leq C_0^*\tau^{1/3} ,\label{dt-err}
\end{align}
where $e^n:= u^n - U^n$.
}
\end{theorem}

The proofs of
Theorem \ref{ErrestFEMSol} and Theorem \ref{ErrestDisSol}  
will be given in Section 3 and Section 4, respectively.
In the rest part of this paper, we denote by $C$ a generic
positive constant which is independent of $\tau$, $h$ and
$n$,
and by $\epsilon$ a generic small positive constant.

\section{Proof of Theorem \ref{ErrestFEMSol}}
\label{secFD}
\setcounter{equation}{0}
In this section, we prove Theorem \ref{ErrestFEMSol}
based on the results of Theorem \ref{ErrestDisSol}.
The proof of the latter is
deferred to Section 4. The following inverse inequalities will be used in this section:
\begin{align}
&\| v \|_{L^p} \le C h^{2/p - 2/q} \| v
\|_{L^q},
~\qquad\quad
\mbox{for}~~~v \in V_h^r , ~~\, 1\leq q\leq p\leq \infty ,
\\[5pt]
 &\| \nabla v \|_{L^p} \le C h^{-1} \| v \|_{L^p},
~\qquad\qquad\mbox{for}~~~v
\in V_h^r ,
~~\, 1\leq p\leq \infty .
\end{align}

\subsection{Preliminaries}

Based on Theorem \ref{ErrestDisSol}, we define
$$
M=\sup_{\tau}\max_{1\leq n\leq N}(\|u^n\|_{W^{1,\infty}}
+\|U^n\|_{W^{1,\infty}})+1
$$
so that
\begin{align*}
&\sigma(|\nabla u^n|^2) \geq \sigma_M ,\,\quad
|\sigma(|\nabla u^n|^2)|+|\sigma'(|\nabla u^n|^2)|
+|\sigma''(|\nabla u^n|^2)|\leq C_M ,\\
&\sigma(|\nabla U^n|^2) \geq \sigma_M ,\quad
|\sigma(|\nabla U^n|^2)|+|\sigma'(|\nabla U^n|^2)|
+|\sigma''(|\nabla U^n|^2)|\leq C_M ,
\end{align*}
for some positive constants $\sigma_M$ and $C_M$.

For any given function $w\in H^1$,
we define the following matrix functions:
\begin{align}
&B(\nabla w) = 2\sigma'(|\nabla w|^2) \nabla w( \nabla
w)^T
\, , \qquad A(\nabla w) = \sigma(|\nabla w|^2) I
+B(\nabla w)
\, . \label{hat-A-n0}
\end{align}
For $n\geq 0$ we define the projection operators $\overline R_h^{n+1}:
H^1(\Omega) \rightarrow V_h^r $ and $R_h^{n+1}: H^1(\Omega)
\rightarrow
V_h^r $ by
\begin{align}
&\big( A(\nabla u^n)\nabla ( w-\overline R_h^{n+1}w),\nabla
v\big)
=0,\quad
\forall~~w\in H^1~~\mbox{and}~~v\in V_h^r , \label{projection0}\\
&\big( A(\nabla U^n)\nabla ( w- R_h^{n+1}w),\nabla v\big)
=0,\quad
\forall~~w\in H^1~~\mbox{and}~~v\in V_h^r ,\label{projection}
\end{align}
where $\int_\Omega \overline R_h^{n+1}w\d x=\int_\Omega
R_h^{n+1}w\d
x=\int_\Omega w\d x$ are enforced for uniqueness, and we set
$\overline R_h^{0}:=\overline R_h^{1}$,
$R_h^{0}:=R_h^{1}$.
These two projection operators are well defined since
\begin{align*}
&\lambda^2\sigma_{M}^3|\xi |^2 \leq \xi^T A(\nabla u^n)\xi \leq
2C_{M}|\xi
|^2, \quad\,\forall~ \xi \in \R^2 \, , \\ 
&\lambda^2\sigma_{M}^3|\xi
|^2 \leq \xi^T A(\nabla U^n)\xi \leq 2C_{M}|\xi |^2,
\quad\forall~
\xi \in \R^2 \, .
\end{align*}

We denote
$$
\theta^{n+1}_h =U^{n+1}-R_h^{n+1}U^{n+1}, \quad \mbox{ and } \quad
\overline\theta^{n+1}_h =u^{n+1}- \overline R_h^{n+1} u^{n+1} 
\, .
$$
By the classical theory of finite element methods, with the
regularity of $U^n$ given in Theorem
\ref{ErrestDisSol},
we have
\begin{align}
&
\|u^{n+1}- R_h^{n+1}u^{n+1}\|_{W^{1,\infty}} \leq
C\|u^{n+1}\|_{H^3}h,\label{proj-001} \\
&\|\overline\theta^{n+1}_h \|_{H^l} \leq
C\|u^{n+1}\|_{H^{r+1}}h^{r+1-l},
\qquad\mbox{for}~~l=0,1,\label{proj-002}  \\
&\| R_h^{n+1}U^{n+1}\|_{W^{1,\infty}}+\| \overline
R_h^{n+1}u^{n+1}\|_{W^{1,\infty}}\leq
C(\|U^{n+1}\|_{W^{1,\infty}}+\|u^{n+1}\|_{W^{1,\infty}}) ,
\label{proj-03}\\
&\|\tau D_\tau\nabla U^{n+1}\|_{L^\infty} \leq C\|\tau
D_\tau\nabla
e^{n+1}\|_{L^\infty}
+C\|\tau D_\tau\nabla u^{n+1}\|_{L^\infty}\leq C\tau^{1/3},
\label{proj-003}\\
&\|D_\tau A(\nabla U^{n})\|_{L^{\bar p}}\leq C\|D_\tau
\nabla
U^{n}\|_{L^{\bar p}}\leq C\|D_\tau
U^{n}\|_{H^2}\label{proj-04} ,
\end{align}
and
\begin{align}
\|\theta^{n+1}_h \|_{W^{l,q}}
&\leq
\|e^{n+1}-R_h^{n+1}e^{n+1}\|_{W^{l,q}}
+\|u^{n+1}-R_h^{n+1}u^{n+1}\|_{W^{l,q}}\nn\\
& \leq Ch^{2-l}
\|e^{n+1}\|_{W^{2,q}}+Ch^{2-l+2/q}
\|u^{n+1}\|_{H^3} \quad\mbox{for $l=0,1$ and $2\leq q\leq p$},
\label{proj-01}
\end{align}
where $p$ is given in Theorem \ref{ErrestDisSol} and $1/\bar p+1/p=1/2$.
The above inequality (\ref{proj-002}) with $l=0,1$ is standard 
$L^2$ and $H^1$ error estimate of the finite element method for
elliptic equations, respectively. 
Since $A(\nabla U^n)\in W^{1,p}$ for some $p>2$, 
the $L^2$ error estimate $\|u^{n+1}- R_h^{n+1}u^{n+1}\|_{L^2} 
\leq Ch^3\|u^{n+1}\|_{H^3}$ is also standard. Then, 
\refe{proj-001} can be derived by introducing an extra interpolation 
and an inverse inequality (see page 93, of the book \cite{Bra}. 
Moreover, \refe{proj-03} and 
(\ref{proj-01}) follow from Theorem 8.1.11 and Theorem 8.5.3 of \cite{BS}, 
respectively, and 
(\ref{proj-003})-(\ref{proj-04}) 
are consequences of Theorem \ref{ErrestDisSol}.  
From these inequalities we also derive that
\begin{align}
\|\theta^{n}_h\|_{W^{1,\infty}} &\leq
\|e^{n}-R_h^{n}e^{n}\|_{W^{1,\infty}}
+\|u^{n}-R_h^{n}u^{n}\|_{W^{1,\infty}}\nn\\
&\leq C
\|e^{n}\|_{W^{1,\infty}}+Ch
\|u^{n}\|_{H^3} \nn\\
&\leq C(\tau^{1/3}+h) .
\label{proj-012}
\end{align}
In this section, we shall frequently use the inequalities
 (\ref{proj-001})-(\ref{proj-012}). Moreover,
we need the following Lemma.

\begin{lemma}\label{LemmP}
{\it Under the assumptions of Theorem {\rm\ref{ErrestFEMSol}},
there exist positive constants $\widetilde\tau_0$ and $\delta_0$ such that when 
$\tau \le \widetilde\tau_0$, 
\begin{align}
&\biggl(\sum_{n=0}^{N-1}\tau\|D_\tau
\theta^{n+1}_h \|_{H^{-1}}^2\biggl)^{\frac{1}{2}}\leq
 C(\tau^{1/3}+h^{\delta_0})h^2, 
\label{proj-06}\\
&\biggl(\sum_{n=0}^{N-1}\tau\|D_\tau
(u^{n+1}-R_h^{n+1}u^{n+1})\|_{H^1}^2\biggl)^{\frac{1}{2}}\leq
Ch^r ,
\label{proj-066}\\
&\biggl(\sum_{n=0}^{N-1}\tau\|D_\tau
\overline\theta^{n+1}_h \|_{L^2}^2\biggl)^{\frac{1}{2}}\leq
Ch^{r+1} .\label{proj-036}
\end{align}
}
\end{lemma}
\noindent{\it Proof}~~ Since $u^n$ is smooth enough,
(\ref{proj-066})-(\ref{proj-036}) can be obtained easily.
Here we only prove  (\ref{proj-06}).
Note that
\begin{align}
&\Big(A(\nabla U^{n})\nabla (U^{n+1}-R_h^{n+1}U^{n+1}), \,
\nabla \phi_h \Big) = 0 ,\\
&\Big(A(\nabla U^{n-1})\nabla (U^{n+1}-R_h^nU^{n+1}), \,
\nabla
\phi_h \Big) = 0 .
\end{align}
The difference of the above two equations gives
\begin{align*}
&\Big(A(\nabla u^{n})\nabla
(R_h^nU^{n+1}-R_h^{n+1}U^{n+1}),
\, \nabla \phi_h \Big) \\
&+\Big((A(\nabla U^{n})-A(\nabla u^{n}))\nabla
(R_h^nU^{n+1}-R_h^{n+1}U^{n+1}),
\, \nabla \phi_h \Big) \\
& + \Big ( (A(\nabla U^{n})-A(\nabla U^{n-1}))\nabla
(U^{n+1}-R_h^nU^{n+1}), \, \nabla \phi_h \Big) =0 ,
\end{align*}
which together with Theorem \ref{ErrestDisSol} implies 
\begin{align*}
&\|\nabla (R_h^nU^{n+1}-R_h^{n+1}U^{n+1})\|_{L^2} \\
&\leq C
\|(A(\nabla
U^{n}) -A(\nabla u^{n}))
\nabla (R_h^nU^{n+1}-R_h^{n+1}U^{n+1})\|_{L^{2}}\\
&~~~+C
\|(A(\nabla
U^{n}) -A(\nabla U^{n-1}))
\nabla (U^{n+1}-R_h^nU^{n+1})\|_{L^{2}}\\
&\leq C\|\nabla e^{n}\|_{L^\infty}
\|\nabla (R_h^nU^{n+1}-R_h^{n+1}U^{n+1})\|_{L^2}
+C \tau\|D_\tau \nabla U^{n}\|_{L^{\bar p}}
\|\nabla (U^{n+1}-R_h^nU^{n+1})\|_{L^p}\\
&\leq C\tau^{1/3}\|\nabla (R_h^nU^{n+1}
-R_h^{n+1}U^{n+1})\|_{L^2}\\
&~~~
+C\tau\|D_\tau U^{n}\|_{H^2}(\|e^{n+1}-R_h^{n}e^{n+1}\|_{W^{1,p}}
+\|u^{n+1}-R_h^{n}u^{n+1}\|_{W^{1,p}})  \\
&\leq C\tau^{1/3}\|\nabla (R_h^nU^{n+1}
-R_h^{n+1}U^{n+1})\|_{L^2}\\
&~~~
+C\tau\|D_\tau U^{n}\|_{H^2}(Ch 
\|e^{n+1}\|_{W^{2,p}}+Ch^{1+2/p}
\|u^{n+1}\|_{H^3})  \\
&\leq C\tau^{1/3}\|\nabla (R_h^nU^{n+1}
-R_h^{n+1}U^{n+1})\|_{L^2}
+C\|D_\tau U^{n}\|_{H^2}(\tau^{1/3}+h^{2/p})\tau h ,
\end{align*}
where we have used (\ref{dt-err}), 
(\ref{proj-04}) and a similar $W^{1,p}$ estimate as given in \refe{proj-01}. 
When $\tau<\widetilde\tau_0:=\min(\tau_0^*,(2C)^{-3})$,
we get
\begin{align}\label{dkm2}
\|\nabla (R_h^nU^{n+1}-R_h^{n+1}U^{n+1})\|_{L^2}
&\leq 2C\|D_\tau U^{n}\|_{H^2}(\tau^{1/3}+h^{2/p})\tau h .
\end{align}

To establish the corresponding $L^2$-norm estimate,
for any given $\varphi\in H^1(\Omega)$ we let $\psi$
be the solution of the equation
$$
-\nabla\cdot\Big(A(\nabla U^{n})\nabla \psi \Big)
=\varphi-\frac{1}{|\Omega|}\int_\Omega\varphi\d x
$$
with the boundary condition $A(\nabla
U^{n})\nabla\psi\cdot{\vec n}=0$ on 
$\partial\Omega$ and $\int_\Omega\psi\d x=0$. Due to the
structure of the matrix $A(\nabla U^n)$, this
boundary condition is equivalent to
$\nabla\psi\cdot{\vec n}=0$ on $\partial\Omega$.
Since $A(\nabla U^{n})$ is uniformly
bounded in $W^{1,p}\cap H^{1+s_0}$, there exists a 
positive constant $\delta_0\in(0,\min(2/p,s_0))$ 
(dependent on the norm $\|\nabla U^n\|_{H^{1+s_0}}$) such that $\|\psi\|_{H^{2+s}}\leq
C\|\varphi\|_{H^s}$ for $s\in[0,\delta_0]$ (see Appendix).  

By noting the fact that $\int_\Omega(R_h^nU^{n+1}-R_h^{n+1}U^{n+1})\d x=0$, we have  
\begin{align*}
&\big(R_h^nU^{n+1}-R_h^{n+1}U^{n+1},\varphi)\\
&=\Big(A(\nabla U^{n})\nabla (R_h^nU^{n+1}
-R_h^{n+1}U^{n+1}), \, \nabla \psi \Big) \\
&=\Big(A(\nabla U^{n}) \nabla
(R_h^nU^{n+1}-R_h^{n+1}U^{n+1}),
\, \nabla (\psi-\Pi_h\psi) \Big)\\
&~~~ - \Big((A(\nabla U^{n})-A(\nabla U^{n-1})) \nabla
(U^{n+1}-R_h^nU^{n+1}),
\, \nabla (\Pi_h\psi-\psi) \Big) \\
&~~~ - \Big( (A(\nabla U^{n})-A(\nabla U^{n-1})
\nabla (U^{n+1}-R_h^nU^{n+1}), \, \nabla \psi \Big): = R_1 + R_2 + R_3   
\end{align*} 

By \refe{proj-01} and \refe{dkm2},   
the first two terms of the right-hand side of the above equation 
are bounded by 
\begin{align*} 
|R_1| & \leq C\|D_\tau U^n\|_{H^2}\|\psi\|_{H^2} 
(\tau^{1/3}+h^{2/p})\tau h^2 ,
\nn \\ 
|R_2| & \leq 
C\|D_\tau A(\nabla  U^n)\|_{L^{\bar p}} 
\|\nabla (U^{n+1}-R_h^nU^{n+1})\|_{L^p}\|\psi\|_{H^2} \tau h  
\nn\\
&\leq C\|D_\tau U^n\|_{H^2}\|\psi\|_{H^2}(\tau^{1/3}+h^{2/p})\tau h^2  ,
\end{align*} 
where $1/\bar p + 1/p=1/2$. 
Again by \refe{dt-err}, \refe{proj-01} and \refe{dkm2} and noting the 
homogeneous boundary condition, with integration by part, 
we can bound the last term by 
\begin{align*}
&|R_3|= \Big|\Big( (A(\nabla U^{n})-A(\nabla U^{n-1}))
\nabla (U^{n+1}-R_h^nU^{n+1}), \, \nabla \psi \Big)\Big|\\
& = \Big|\Big( 
U^{n+1}-R_h^nU^{n+1} , 
\nabla\cdot\big[(A(\nabla U^{n})-A(\nabla U^{n-1}))\nabla \psi\big]    \Big)   \Big|
\\
&\leq  \|U^{n+1}-R_h^nU^{n+1}\|_{L^p} 
\|\nabla\cdot\big[(A(\nabla U^{n})-A(\nabla U^{n-1}))\nabla \psi\big] \|_{L^{p'}} 
\\ 
&\leq C(h^2\|e^{n+1}\|_{W^{2,p}}+h^{2+2/p}\|u^{n+1}\|_{H^3})
\left ( 
\|A(\nabla U^{n})-A(\nabla U^{n-1})\|_{H^1} \| \nabla\psi\|_{L^{\widetilde p}} 
\right. \nn \\ 
&~~~+ \left. 
\|A(\nabla U^{n})-A(\nabla U^{n-1})\|_{L^{\widetilde p}}\|\psi\|_{H^2} \right ) 
\nn \\
&\leq C\|D_\tau U^n\|_{H^2}\|\varphi\|_{L^2}
(\tau^{1/3}+h^{2/p})\tau h^2 ,
\end{align*}
where $1/p+1/p'=1 $ and $1/\widetilde p + 1/2=1/p' $.

With the above estimates, we obtain 
\begin{align*}
&\|R_h^nU^{n+1}-R_h^{n+1}U^{n+1}\|_{L^2} \leq C\|D_\tau
U^n\|_{H^2}(\tau^{1/3}+h^{2/p}) \tau h^2
\quad\mbox{for}~~n\geq 1.
\end{align*}
Since $R_h^0U^1=R_h^1U^1$, we have   
\begin{align*}
\biggl(\sum_{n=0}^{N-1}\tau\|R_h^nU^{n+1}
-R_h^{n+1}U^{n+1}\|_{L^2}^2\biggl)^{\frac{1}{2}}\leq
C(\tau^{1/3}+h^{2/p})\tau h^2 .
\end{align*}

Finally, we take a standard approach to the $H^{-1}$-norm estimate 
\refe{proj-06} \cite{BS}. Since 
\begin{align*}
|(\phi-R_h\phi,\varphi)|
&=\inf_{\psi_h\in V_h^r}|(A(\nabla U^n)\nabla(\phi-R_h\phi),\nabla(\psi-\psi_h))|\\
&\leq C\|\nabla(\phi-R_h\phi)\|_{L^2}\|\psi\|_{H^{2+\delta_0}}h^{1+\delta_0} \\
&\leq C\|\phi\|_{H^2} 
\|\varphi\|_{H^{\delta_0}}h^{2+\delta_0} ,\quad\forall~\varphi\in H^{\delta_0} ,
\end{align*}
we have 
\begin{align} 
\|\phi-R_h\phi\|_{H^{-\delta_0}}\leq C\|\phi\|_{H^2}h^{2+\delta_0} ,
\quad \forall~\phi\in H^2 
\nn 
\end{align}
from which, we further derive that
\begin{align}
&\|D_\tau (U^{n+1} -R_h^{n+1}U^{n+1})\|_{H^{-1}}  
\nn \\
&\leq \|D_\tau U^{n+1} -R_h^nD_\tau U^{n+1}\|_{H^{-1}}
+\tau^{-1} \|R_h^{n+1}U^{n+1}-R_h^nU^{n+1}\|_{H^{-1}} \nn\\
&\leq \|D_\tau U^{n+1} -R_h^nD_\tau U^{n+1}\|_{H^{-\delta_0}}
+\tau^{-1}
\|R_h^{n+1}U^{n+1}-R_h^nU^{n+1}\|_{L^2} \nn\\
&\leq C\|D_\tau e^{n+1}\|_{H^2}h^{2+\delta_0}+C\|D_\tau
u^{n+1}\|_{H^3}h^3 +
\tau^{-1}
\|R_h^{n+1}U^{n+1}-R_h^nU^{n+1}\|_{L^2} . \nn
\end{align}
(\ref{proj-06}) follows immediately.
~\endproof
\medskip



\subsection{Boundedness of the numerical solution}
\label{SecBDN}
By \refe{dt-bound} and \refe{proj-03}, we can re-define
\begin{align*}
M&=\sup_{\tau,h}\Big(\max_{0\leq n\leq
N}\|u^n\|_{W^{1,\infty}}
+\max_{0\leq n\leq
N}\| \overline R_h^n u^n\|_{W^{1,\infty}}\nn\\
&~~~+\max_{0\leq n\leq N}\|U^n\|_{W^{1,\infty}} +\max_{0\leq
n\leq
N}\| R_h^n U^n\|_{W^{1,\infty}}\Big)+2 .
\end{align*}
By the regularity assumptions on $\sigma$, there
exist $\sigma_{M}$ and $C_{M}>0$ such that
\begin{align}
&\sigma(s^2)\geq \sigma_{M} ,
\qquad\qquad\qquad\qquad\qquad\forall~s\in[-M,M], \label{K122}\\[5pt]
&|\sigma(s^2)|+|\sigma'(s^2)|+|\sigma''(s^2)| \leq C_{M},
\quad~\forall~s\in[-M,M] . \label{K133}
\end{align}

\begin{lemma}
\label{l-3-3} {\it
Under the assumptions of Theorem {\rm\ref{ErrestFEMSol}}, 
there exist positive constants $\widehat\tau_0$ and $\widehat h_0$
which are independent of $n$, $\tau$ and $h$, such that the finite element system
{\rm(\ref{FDFEM})}
admits a unique solution $\{ U^n_h \}_{n=1}^N$ when $
\tau<\widehat\tau_0$ and $h<\widehat h_0$, satisfying
\begin{align}
&\|U_h^n\|_{L^\infty}+\|\nabla U_h^n\|_{L^\infty}\leq M
,\label{Uh-bound} \\
&\|e^n_h\|_{L^\infty}+\|\nabla e^n_h\|_{L^\infty}
<\tau^{1/8}+h^{\delta_0/8} \label{Uh-bound2} ,
\end{align}
where $e_h^{n}=R_h^{n}U^{n}-U^{n}_h$ and 
$\delta_0$ is given in Lemma {\rm\ref{LemmP}}.
}
\end{lemma}

\noindent{\it Proof}~~~
By (\ref{K122})-(\ref{K133}), the coefficient matrix of the linear
system (\ref{FDFEM}) is symmetric and positive definite, which implies that
(\ref{FDFEM}) admits a unique solution
$U^{n+1}_h\in V_h^r $ for $0\leq n\leq k$.

It is easy to see that the inequalities
\refe{Uh-bound}-\refe{Uh-bound2} hold for $n=0$.  By  mathematical induction,
we can assume that \refe{Uh-bound}-(\ref{Uh-bound2}) hold for $0\leq n\leq
k$ for some $k\geq 0$.

Since the solution $U^{n+1}$ of (\ref{e-TDparab-1}) satisfies
\begin{align}
&\big(D_\tau U^{n+1},v\big)+\big(\sigma(|\nabla U^n|^2)
\nabla
U^{n+1},\nabla v\big) =\big(g^{n+1},v\big) ,\quad
\forall~v\in V_h^r
, \nn
\end{align}
the error function $e_h^{n+1}$ satisfies
\begin{align}
&\big(D_\tau e^{n+1}_h,v\big)+\big(\sigma(|\nabla U^n|^2)
\nabla
e^{n+1}_h,\nabla v\big)
\label{errFEMFDFEM}\\
& = \Big[ -\big( \sigma(|\nabla U^n|^2)
\nabla \theta^{n+1}_h , \, \nabla v \big) +\big((\sigma(|\nabla U_h^n|^2)
- \sigma(|\nabla U^n|^2))\nabla
U^{n+1}_h,\nabla v\big) \Big]  -\big(D_\tau \theta^{n+1}_h ,v\big)
\nn \\
& := J_1(v) + J_2(v), \qquad\forall~ v\in V_h^r . \nn
\end{align}
By using Taylor's expansion, we see that
\begin{align}
& (\sigma(|\nabla U_h^n|^2) - \sigma(|\nabla U^n|^2))
 \nabla U^{n+1}_h \nn\\
&= \left ( 2\sigma'(|\nabla U^n|^2)\nabla U^n \cdot \nabla
(U_h^n-U^n) +\sigma'(|\nabla U^n|^2)|\nabla(U_h^n-U^n)|^2
\right )
\nabla U^{n+1}
\nn \\
&~~~+\frac{1}{2}\sigma''(\xi^n_h) |\nabla (U_h^n+
U^n)\cdot\nabla
(U_h^n-U^n)|^2 \nabla U^{n+1} \nn\\
&~~~+(\sigma(|\nabla U_h^n|^2) - \sigma(|\nabla U^n|^2))
\nabla
 (-e^{n+1}_h-\theta^{n+1}_h )\nn\\
& =-2 \sigma'(|\nabla U^n|^2) \nabla U^n \cdot\nabla(
e^n_h+\theta^{n+1}_h ) (\nabla U^n+\tau
D_\tau\nabla U^{n+1}) \nn\\
&~~~+2\sigma'(|\nabla U^n|^2) \nabla U^n \cdot\nabla\tau
D_\tau\theta^{n+1}_h  \nabla U^{n+1} \nn\\
&~~~+ \Big(\sigma'(|\nabla
U^n|^2)|\nabla(U_h^n-U^n)|^2+\frac{1}{2}\sigma''(\xi^n_h)
|\nabla(
U^n+U_h^n)\cdot\nabla ( e^n_h+\theta^n_h )|^2 \Big)
\nabla U^{n+1} 
\nn \\
&~~~ - (\sigma(|\nabla U_h^n|^2) - \sigma(|\nabla U^n|^2))
\nabla ( e^{n+1}_h+\theta^{n+1}_h ) \label{sigma-h-1}
\end{align}
where $\xi_h^n$ is some number between $|\nabla U_h^n|^2$ and
$|\nabla U^n|^2$. By using the notations in (\ref{hat-A-n0}), we see further that
\begin{align*}
J_1(v)
&= -\big( A(\nabla U^n) \nabla \theta^{n+1}_h ,
\, \nabla v \big)
\nn \\
&~~~-\big(2 \sigma'(|\nabla U^n|^2) \nabla U^n \cdot\nabla
\theta^{n+1}_h \tau D_\tau\nabla U^{n+1} ,
\, \nabla v \big)\nn\\
&~~~-\big(2 \sigma'(|\nabla U^n|^2) (\nabla U^n \cdot\nabla
e^n_h)
\nabla U^{n+1},
\, \nabla v \big)\nn\\
&~~~+\big(2 \sigma'(|\nabla U^n|^2) \nabla U^n \cdot\nabla\tau
D_\tau \theta^{n+1}_h  \nabla U^{n+1}, \, \nabla v \big)
\nn\\
&~~~+ \big(\sigma'(|\nabla
U^n|^2)|\nabla(U_h^n-U^n)|^2+\frac{1}{2}\sigma''(\xi^n_h) |\nabla(
U^n+U_h^n)\cdot\nabla (
e^n_h+\theta^n_h )|^2 \nabla
U^{n+1}, \, \nabla v \big) 
\nn \\
&~~~ - \big((\sigma(|\nabla U_h^n|^2) - \sigma(|\nabla
U^n|^2))
\nabla ( e^{n+1}_h+\theta^{n+1}_h ) , \, \nabla v
\big).
\end{align*}
Let
\begin{align}
\gamma(|\nabla U^n|^2)=2|\sigma'(|\nabla U^n|^2)|
|\nabla U^n|^2=\sigma(|\nabla U^n|^2)-\lambda^2\sigma(|\nabla U^n|^2)^3.
\label{gamma}
\end{align}
Taking $v=e_h^{n+1}$ in (\ref{errFEMFDFEM}) and noting the fact 
$\big( A(\nabla U^n) \nabla \theta^{n+1}_h ,
\, \nabla e_h^{n+1} \big)=0$, we obtain
\begin{align*}
J_1(e_h^{n+1})
& \leq \big(\gamma(|\nabla U^n|^2)|
\nabla e_h^{n+1}|,|\nabla
e_h^{n}|\big) \nn\\
&~~~ +C\|\tau D_\tau\nabla U^{n+1}\|_{L^\infty}(\|\nabla
e_h^{n}
\|_{L^2}+\| \nabla \theta^{n+1}_h
\|_{L^2})\|\nabla e_h^{n+1} \|_{L^2} \\
 &~~~ + C\biggl(\sum_{m=n}^{n+1}\|e^{m}- R_h^{m}
e^{m}\|_{H^1}+ \tau \|D_\tau(u^{n+1}- R_h^{n+1} u^{n+1})\|_{H^1}\biggl)  
\|\nabla e_h^{n+1} \|_{L^2} \\
&~~~+C(\|\nabla e_h^{n} \|_{L^\infty}+\|\nabla \theta^{n}_h
\|_{L^\infty})(\|\nabla e_h^{n} \|_{L^2}+\|\nabla \theta^{n}_h
\|_{L^2})\|\nabla e_h^{n+1} \|_{L^2} \\
&~~~+C(\|\nabla e_h^{n} \|_{L^\infty}+\|\nabla \theta^{n}_h
\|_{L^\infty})(\|\nabla \theta^{n+1}_h \|_{L^2}\|\nabla e_h^{n+1}
\|_{L^2}+\|\nabla e_h^{n+1} \|_{L^2}^2) .
\end{align*}
From \refe{proj-01}, \refe{proj-012} and \refe{Uh-bound2} we have 
\begin{align*}
&\|\nabla \theta^{n}_h\|_{L^2}\leq Ch\|e^n\|_{H^2}+Ch^2,\\
&\|\nabla e_h^{n} \|_{L^\infty}+\|\nabla \theta^{n}_h
\|_{L^\infty}\leq C(\tau^{1/8}+h^{\delta_0/8})<\epsilon
\end{align*}
when $\tau<\tau_1$ and $h<h_1$ for some positive constants $\tau_1$ and $h_1$ 
(which depend on the constant $\epsilon$). With 
(\ref{proj-001})-(\ref{proj-012}), the induction assumptions
(\ref{Uh-bound})-(\ref{Uh-bound2}) and the regularity of
$U^n$ given
in Theorem \ref{ErrestDisSol}, we derive that, 

\begin{align*}
J_1(e_h^{n+1})& \leq \frac{1}{2}\Big\|\sqrt{\gamma(|\nabla U^n|^2)}
\nabla e_h^{n+1}\Big \|_{L^2}^2
+\frac{1}{2}\Big\|\sqrt{\gamma(|\nabla U^n|^2)} \nabla e_h^{n} \Big\|_{L^2}^2
+C\tau^{1/3}\|\nabla e_h^{n} \|_{L^2}\|\nabla e_h^{n+1} \|_{L^2}\\
&~~~ + C(h\|e^{n+1}\|_{H^2}+h\|e^{n}\|_{H^2}+\tau^{1/3}
h^2+\tau
\|D_\tau(u^{n+1}- R_h^{n+1} u^{n+1})\|_{H^1})\|\nabla e_h^{n+1} \|_{L^2}
 \\
&~~~+\epsilon(\|\nabla e_h^{n} \|_{L^2}^2
+\|\nabla e_h^{n+1}\|_{L^2}^2) 
+C\epsilon^{-1}(\|\nabla e_h^{n} \|_{L^\infty}^2+\|\nabla
\theta^{n}_h
\|_{L^\infty}^2)(\|\nabla
\theta^{n}_h\|_{L^2}^2+\|\nabla
\theta^{n+1}_h \|_{L^2}^2)\\
& \leq \frac{1}{2}\Big\|\sqrt{\gamma(|\nabla U^n|^2)}
\nabla e_h^{n+1}\Big \|_{L^2}^2
+\frac{1}{2}\Big\|\sqrt{\gamma(|\nabla
U^n|^2)} \nabla e_h^{n} \Big\|_{L^2}^2\\
&~~~+C\epsilon^{-1}(h^2\|e^{n+1}\|_{H^2}^2+h^2\|e^{n}\|_{H^2}^2+\tau^{2/3}
h^4+\tau^2
\|D_\tau(u^{n+1}- R_h^{n+1} u^{n+1})\|_{H^1}^2)\\
&~~~+2\epsilon(\|\nabla e_h^{n} \|_{L^2}^2
+\|\nabla e_h^{n+1} \|_{L^2}^2) \\
&~~~
+C\epsilon^{-1}(\|\nabla e_h^{n} \|_{L^\infty}^2
+\tau^{2/3}+h^2)(h^2\|e^{n}\|_{H^2}^2+h^2\|e^{n+1}\|_{H^2}^2+h^4)\\
& \leq \frac{1}{2}\Big\|\sqrt{\gamma(|\nabla U^n|^2)}
\nabla e_h^{n+1}\Big \|_{L^2}^2
+\frac{1}{2}\Big\|\sqrt{\gamma(|\nabla
U^n|^2)} \nabla e_h^{n} \Big\|_{L^2}^2\\
&~~~+3\epsilon(\|\nabla e_h^{n} \|_{L^2}^2
+\|\nabla e_h^{n+1} \|_{L^2}^2)+C\epsilon^{-1}(h^2\|e^{n}\|_{H^2}^2
+h^2\|e^{n+1}\|_{H^2}^2)\\
&~~~+C\epsilon^{-1}(\tau^{2/3}h^4+h^6)+C\epsilon^{-1}\tau^2
\|D_\tau(u^{n+1}- R_h^{n+1} u^{n+1})\|_{H^1}^2 ,
\end{align*}
where we have used the inverse inequality 
$h^4\|\nabla e^n_h\|_{L^\infty}^2
\leq Ch^2\|\nabla e^n_h\|_{L^2}^2
\leq \epsilon \|\nabla e^n_h\|_{L^2}^2$. 
For $J_2(e_h^{n+1})$, we have the following estimate,  
\begin{align*}
J_2(e_h^{n+1})
& \leq \| D_\tau\theta^{n+1}_h
\|_{H^{-1}}\|e_h^{n+1} \|_{H^1} \\
&\leq C\epsilon^{-1}\| D_\tau\theta^{n+1}_h
\|_{H^{-1}}^2+\epsilon\|\nabla e_h^{n+1}
\|_{L^2}^2+\epsilon\|e_h^{n+1} \|_{L^2}^2 . 
\end{align*}
With the above estimates, \refe{errFEMFDFEM} reduces to
\begin{align}
&\frac{1}{2} D_\tau \| e^{n+1}_h \|_{L^2}^2 + \frac{1}{2}
\Big(\big\|\sqrt{\sigma(|\nabla U^n|^2)}\nabla
e^{n+1}_h\big\|_{L^2}^2-\big\|\sqrt{\gamma(|\nabla
U^n|^2)}\nabla e^{n}_h\big\|_{L^2}^2\Big) \nn\\
& \leq 3\epsilon (\| \nabla e^n
\|_{L^2}^2+\|
\nabla e^{n+1} \|_{L^2}^2)
+C\| D_\tau \theta^{n+1}_h \|_{H^{-1}}^2 \nn\\
&~~~+C\epsilon^{-1}\|e_h^{n+1}\|_{L^2}^2+ C\epsilon^{-1}\tau^2
\|D_\tau(u^{n+1}- R_h^{n+1}u^{n+1})\|_{H^1}^2  \nn\\
&~~~+
C\epsilon^{-1}(\|e^{n}\|_{H^2}^2+\|e^{n+1}\|_{H^2}^2)h^{2}
+C\epsilon^{-1}(\tau^{2/3}h^4+h^{6}) . \label{dkhf2}
\end{align}
From (\ref{dt-err}) we derive that
$$
\|\tau D_\tau\gamma(|\nabla U^n|^2)\|_{L^\infty}
\leq C \|\tau D_\tau e^n\|_{W^{1,\infty}}
+C\|\tau D_\tau u^n\|_{W^{1,\infty}} \leq C\tau^{1/3} ,
$$
which implies 
\begin{align*}
& \| \sqrt{\sigma(|\nabla U^n|^2)} \nabla e_h^{n+1} \|_{L^2}^2
- \|
\sqrt{\gamma(|\nabla U^n|^2)} \nabla e_h^{n} \|_{L^2}^2
\\
&=\| \sqrt{\sigma(|\nabla U^n|^2)-\gamma(|\nabla U^n|^2)}
\nabla
e_h^{n+1} \|_{L^2}^2 +\| \sqrt{\gamma(|\nabla U^{n}|^2)}
\nabla
e_h^{n+1} \|_{L^2}^2- \| \sqrt{\gamma(|\nabla U^{n-1}|^2)}
\nabla
e_h^{n} \|_{L^2}^2
\nn\\
&~~~ - \big((\gamma(|\nabla U^n|^2)-\gamma(|\nabla
U^{n-1}|^2)) \nabla e_h^{n},\nabla e_h^{n}\big) \\
& \geq \|\lambda \sigma(|\nabla U^n|^2)^{3/2}
\nabla
e_h^{n+1} \|_{L^2}^2+\tau D_\tau\| \sqrt{\gamma(|\nabla
U^n|^2)}
\nabla e_h^{n+1} \|_{L^2}^2 - \tau \| \sqrt{|D_{\tau}
\gamma(|\nabla
U^n|^2)|} \nabla e_h^{n} \|_{L^2}^2 \\
& \geq \lambda^2\sigma_M^3\|\nabla
e_h^{n+1} \|_{L^2}^2+\tau D_\tau\| \sqrt{\gamma(|\nabla
U^n|^2)}
\nabla e_h^{n+1} \|_{L^2}^2 - C\tau^{1/3} \|\nabla e_h^{n} \|_{L^2}^2.
\end{align*}
With the above inequality, (\ref{dkhf2}) reduces to
\begin{align*}
&\frac{1}{2} D_\tau \| e^{n+1}_h \|_{L^2}^2
+ \frac{\lambda^2\sigma_M^3}{2}\|\nabla
e_h^{n+1} \|_{L^2}^2+ \frac{\tau}{2} D_\tau\| \sqrt{\gamma(|\nabla
U^n|^2)}
\nabla e_h^{n+1} \|_{L^2}^2 \\
& \leq 3\epsilon  (\| \nabla e^n
\|_{L^2}^2+\|
\nabla e^{n+1} \|_{L^2}^2)
+C\| D_\tau \theta^{n+1}_h \|_{H^{-1}}^2\\
&~~~+C\epsilon^{-1}\|e_h^{n+1}\|_{L^2}^2+
C\epsilon^{-1}\tau^2
\|D_\tau(u^{n+1}- R_h^{n+1}u^{n+1})\|_{H^1}^2 \\
&~~~+
C\epsilon^{-1}(\|e^{n}\|_{H^2}^2+\|e^{n+1}\|_{H^2}^2)h^{2}
+C\epsilon^{-1}(\tau^{2/3}h^4+h^{6}) .
\end{align*} 

Choosing $\epsilon =\lambda^2\sigma_M^3/72$, 
by Theorem \ref{ErrestDisSol},  Lemma \ref{LemmP}
and Gronwall's inequality, we derive that
\begin{align*}
&\|e^{k+1}_h\|_{L^2}^2 + \sum_{m=0}^k\tau \|\nabla
e^{m+1}_h\|_{L^2}^2
\leq
C\tau^{2/3}h^4+Ch^{4+2\delta_0} .
\end{align*}
when $\tau<\tau_2 \le \widetilde \tau_0$ 
and $h<h_2$ for some positive constants $\tau_2$ and $h_2$.

By the inverse inequality, we have
\begin{align*}
\|e^{k+1}_h\|_{L^\infty}+\|\nabla e^{k+1}_h\|_{L^\infty}&
\leq C h^{-2} \|e^{k+1}_h\|_{L^2}
\leq C(\tau^{1/3}+h^{\delta_0}) \, ,
\end{align*}
which implies that
\begin{align}
\|e^{k+1}_h\|_{L^\infty}+\|\nabla e^{k+1}_h\|_{L^\infty}
\leq \tau
^{1/8}+h^{\delta_0/8} ,
\end{align}
and 
\begin{align}
\|U^{k+1}_h\|_{L^\infty}+\|\nabla U^{k+1}_h\|_{L^\infty}
\leq
\|R^{k+1}_hU^{k+1}\|_{L^\infty}+\|\nabla
R^{k+1}_hU^{k+1}\|_{L^\infty}+1\leq M 
\end{align}
when $\tau<\tau_3$ and $h<h_3$ for some positive constants $\tau_3$ and $h_3$.
The induction on
(\ref{Uh-bound})-(\ref{Uh-bound2}) is closed with 
$\widehat\tau_0=\min\{\tau_0^*,\widetilde\tau_0,\tau_1,\tau_2,\tau_3\}$ 
and $\widehat h_0=\min\{h_1,h_2,h_3\}$. 

The proof of Lemma \ref{l-3-3} is completed. ~\endproof

\subsection{Unconditionally optimal error estimate}
Now we turn back to the proof of Theorem \ref{ErrestFEMSol}.
Let $\overline e_h^{n}=\overline R_h^{n}u^{n}-U^{n}_h$.
From Lemma \ref{l-3-3}, Theorem \ref{ErrestDisSol},
(\ref{proj-001}) and (\ref{proj-012}), we see that there exist
positive constants $\tau_4<\widehat\tau_0$ and $h_4<\widehat h_0$ such that when $\tau<\tau_4$
and $h<h_4$ 
\begin{align}
&\|U_h^n\|_{L^\infty}+\|\nabla U_h^n\|_{L^\infty}\leq M
,\qquad\qquad\qquad~\mbox{for $n=0,1,\cdots,N$}, \label{Uhbd01} \\
&\|\overline e^n_h\|_{L^\infty}+\|\nabla\overline
e^n_h\|_{L^\infty}
<2\tau^{1/8}+2h^{\delta_0/8}\label{Uhbd02} ,\qquad\mbox{for $n=0,1,\cdots,N$}.
\end{align}

Since the exact solution $u^n$ satisfies
\begin{align}
&\big(D_\tau u^{n+1},v\big)+\big(\sigma(|\nabla u^n|^2)
\nabla
u^{n+1},\nabla v\big) =\big(g^{n+1},v\big)+\big({\cal
E}^{n+1}_{\rm
tr},v\big) ,\quad \forall~v\in V_h^r  , \nn
\end{align}
the error function $\overline{e}_h^{n+1}$ satisfies
\begin{align}
&\big(D_\tau\overline e^{n+1}_h,v\big)+\big(\sigma(|\nabla u^n|^2)
\nabla\overline  e^{n+1}_h,\nabla v\big)
\label{errFE2}\\
& = \Big[ -\big( \sigma(|\nabla u^n|^2)
 \nabla\overline\theta^{n+1}_h , \, \nabla v \big)
+\big((\sigma(|\nabla U_h^n|^2) - \sigma(|\nabla u^n|^2))\nabla U^{n+1}_h,\nabla v\big)
\Big]\nn\\
&~~~
- \big(D_\tau \overline\theta^{n+1}_h ,v\big)
+\big({\cal E}^{n+1}_{\rm tr},v\big)
\nn \\
& := \overline J_1(v) +\overline  J_2(v) + \overline J_3(v),
\qquad\forall~ v\in V_h^r . \nn
\end{align}
To estimate $\overline J_i$, $i=1,2,3,$ we take the same approach
as used for $J_1$ and $J_2$ in Section \ref{SecBDN} and we get
\begin{align*}
 \overline J_1(\overline e_h^{n+1})
  &=-\big( A(\nabla u^n)
 \nabla \overline\theta^{n+1}_h ,
\, \nabla \overline e_h^{n+1} \big)
\nn \\
&~~~-\big(2 \sigma'(|\nabla u^n|^2) \nabla u^n \cdot\nabla
\overline
e^n_h \nabla u^{n+1},
\, \nabla \overline e_h^{n+1} \big)\nn\\
&~~~-\big(2 \sigma'(|\nabla u^n|^2) \nabla u^n \cdot\nabla
\overline\theta^{n+1}_h  \tau D_\tau\nabla u^{n+1} ,
\, \nabla \overline e_h^{n+1} \big)\nn\\
&~~~+\big(2 \sigma'(|\nabla u^n|^2) \nabla u^n \cdot\nabla
\tau
D_\tau\overline\theta^{n+1}_h  \nabla u^{n+1},
\, \nabla \overline e_h^{n+1} \big)
\nn\\
&~~~+ \biggl(\sigma'(|\nabla u^n|^2)|\nabla(u^n_h-u^n)|^2
+\frac{1}{2}\sigma''(\overline \xi^n_h)
|\nabla(u^n+u^n_h)\cdot\nabla (\overline
e^n_h+\overline\theta^{n}_h)|^2 \nabla u^{n+1}, \, \nabla \overline e_h^{n+1}
\biggl)
\nn \\
&~~~ - \big((\sigma(|\nabla U_h^n|^2) - \sigma(|\nabla
u^n|^2))
\nabla ( \overline e^{n+1}_h+\overline\theta^{n+1}_h ),
 \, \nabla \overline e_h^{n+1} \big) \nn\\
& \leq \frac{1}{2}\Big\|\sqrt{\gamma(|\nabla u^n|^2)}
\nabla\overline e_h^{n+1}\Big \|_{L^2}^2
+\frac{1}{2}\Big\|\sqrt{\gamma(|\nabla u^n|^2)} \nabla
\overline e_h^{n} \Big\|_{L^2}^2+C\tau\|\nabla \overline e_h^{n}
\|_{L^2}\|\nabla \overline e_h^{n+1} \|_{L^2} \\
&~~~ + (C\tau+C\tau \|D_\tau\overline\theta^{n+1}_h \|_{H^1})
\|\nabla \overline e_h^{n+1}\|_{L^2} \\
&~~~+C(\|\nabla\overline  e_h^{n} \|_{L^\infty}+\|\nabla
\overline\theta^{n}_h \|_{L^\infty})(\|\nabla \overline
e_h^{n}\|_{L^2}+\|\nabla \overline\theta^{n}_h
\|_{L^2})\|\nabla\overline e_h^{n+1} \|_{L^2}\\
&~~~+C(\|\nabla\overline  e_h^{n} \|_{L^\infty}+\|\nabla
(u^{n}-\overline R_h^nu^n) \|_{L^\infty})\|\nabla\overline
e_h^{n+1} \|_{L^2}^2\\
&~~~+C\|\nabla \overline\theta^{n+1}_h
\|_{L^\infty}(\|\nabla\overline  e_h^{n} \|_{L^2}+\|\nabla
\overline\theta^{n}_h\|_{L^2})\|\nabla\overline
e_h^{n+1}\|_{L^2}\\
& \leq \frac{1}{2}\Big\|\sqrt{\gamma(|\nabla u^n|^2)}
\nabla\overline e_h^{n+1}\Big \|_{L^2}^2
+\frac{1}{2}\Big\|\sqrt{\gamma(|\nabla
u^n|^2)} \nabla\overline e_h^{n} \Big\|_{L^2}^2 \\
&~~~+\epsilon (\|\nabla\overline
e_h^{n}\|_{L^2}^2 +\|\nabla\overline e_h^{n+1}\|_{L^2}^2) +C\epsilon^{-1}(1+
\|D_\tau\overline\theta^{n+1}_h \|_{H^1}^2) \tau^2
+C\epsilon^{-1}h^{2r+2} , 
\end{align*}
\begin{align*}
\overline J_2(\overline e_h^{n+1}) &\leq C\epsilon^{-1}
\| D_\tau \overline\theta^{n+1}_h \|_{L^2}^2
+\epsilon\|\overline e_h^{n+1} \|_{L^2}^2 ,
\end{align*}
and
\begin{align*}
\overline J_3(\overline e_h^{n+1}) &\leq
\epsilon\|\overline e_h^{n+1} \|_{L^2}^2
+C\epsilon^{-1}\|{\cal E}_{\rm tr}^{n+1}\|_{L^2}^2 
\end{align*}
when $\tau<\tau_5$ and $h<h_5$ for some positive constants $\tau_5$ and $h_5$.
With the above estimates, \refe{errFE2} reduces to
\begin{align}
&\frac{1}{2} D_\tau \|\overline e^{n+1}_h \|_{L^2}^2 +
\frac{1}{2}
\Big(\big\|\sqrt{\sigma(|\nabla U^n|^2)}\nabla \overline
e^{n+1}_h\big\|_{L^2}^2-\big\|\sqrt{\gamma(|\nabla
U^n|^2)}\nabla\overline  e^{n}_h\big\|_{L^2}^2\Big) \label{lkh3}\\
& \leq \epsilon (\| \nabla\overline
e^n\|_{L^2}^2+\|\nabla\overline  e^{n+1} \|_{L^2}^2)
+\epsilon\|\overline e_h^{n+1}\|_{L^2}^2\nn\\
&~~~+C\epsilon^{-1}\tau^2 \|D_\tau\overline\theta^{n+1}_h \|_{H^1}^2
+C\epsilon^{-1}\| D_\tau \overline\theta^{n+1}_h \|_{L^2}^2
+C\epsilon^{-1}\|{\cal E}_{\rm tr}^{n+1}\|_{L^2}^2 +
C\epsilon^{-1}(\tau^2+h^{2r+2})  . \nn
\end{align}
Since
\begin{align*}
& \| \sqrt{\sigma(|\nabla U^n|^2)} \nabla \overline e_h^{n+1} \|_{L^2}^2
- \|
\sqrt{\gamma(|\nabla U^n|^2)} \nabla \overline e_h^{n} \|_{L^2}^2
\\
&=\| \sqrt{\sigma(|\nabla U^n|^2)-\gamma(|\nabla U^n|^2)}
\nabla
\overline e_h^{n+1} \|_{L^2}^2 +\| \sqrt{\gamma(|\nabla U^{n}|^2)}
\nabla
\overline e_h^{n+1} \|_{L^2}^2- \| \sqrt{\gamma(|\nabla U^{n-1}|^2)}
\nabla
\overline e_h^{n} \|_{L^2}^2
\nn\\
&~~~ - \big((\gamma(|\nabla U^n|^2)-\gamma(|\nabla
U^{n-1}|^2)) \nabla \overline e_h^{n},\nabla \overline e_h^{n}\big) \\
& \geq \|\lambda\sigma(|\nabla U^n|^2)^{3/2}
\nabla\overline
e_h^{n+1} \|_{L^2}^2+\tau D_\tau\| \sqrt{\gamma(|\nabla
U^n|^2)}
\nabla \overline e_h^{n+1} \|_{L^2}^2 - \tau \| \sqrt{|D_{\tau}
\gamma(|\nabla
U^n|^2)|} \nabla \overline e_h^{n} \|_{L^2}^2 \\
& \geq \lambda^2\sigma_M^3\|\nabla\overline
e_h^{n+1} \|_{L^2}^2+\tau D_\tau\| \sqrt{\gamma(|\nabla
U^n|^2)}
\nabla\overline  e_h^{n+1} \|_{L^2}^2 - C\tau^{1/3}
\|\nabla \overline e_h^{n} \|_{L^2}^2,
\end{align*}
the inequality (\ref{lkh3}) reduces to
\begin{align}
&\frac{1}{2} D_\tau \|\overline e^{n+1}_h \|_{L^2}^2 +
\frac{\lambda^2\sigma_M^3}{2}
\|\nabla \overline
e^{n+1}_h\|_{L^2}^2+ \frac{\tau}{2} D_\tau\| \sqrt{\gamma(|\nabla
U^n|^2)}
\nabla e_h^{n+1} \|_{L^2}^2 \\
& \leq \epsilon (\| \nabla\overline e^n\|_{L^2}^2+\|
\nabla\overline  e^{n+1} \|_{L^2}^2)
+\epsilon\|\overline e_h^{n+1}\|_{L^2}^2 \nn\\
&~~~+C\epsilon^{-1}\tau^2
\|D_\tau\overline\theta^{n+1}_h \|_{H^1}^2
+C\epsilon^{-1}\| D_\tau \overline\theta^{n+1}_h \|_{L^2}^2
+C\epsilon^{-1}\|{\cal E}_{\rm tr}^{n+1}\|_{L^2}^2 +
C\epsilon^{-1}(\tau^2+h^{2r+2}) \, .
\nn
\end{align}
By choosing $\epsilon =\lambda^2\sigma_M^3/24$ and applying Gronwall's inequality, when
$\tau<\tau_6$ and $h <h_6$ for some positive constants $\tau_6$ and $h_6$, we obtain
\begin{align}
\label{error-l2} \max_{0\leq n\leq N}\|\overline e^{n}_h\|_{L^2}^2 +
\sum_{n=0}^N\tau
\|\nabla \overline e^{n}_h\|_{L^2}^2 &\leq
C(\tau^2+h^{2r+2}) .
\end{align}

So far we have proved Theorem \ref{ErrestFEMSol} for the
case $\tau<\tau_7:=\min\{\tau_4,\tau_5,\tau_6\}$ and $h<h_7:=\min\{h_4,h_5,h_6\}$.
Now we consider the case that $\tau\geq\tau_7$ or $h\geq h_7$.
Substituting $v=U^{n+1}_h$ in (\ref{FDFEM}), we get
\begin{align*}
&D_\tau\biggl(\frac{1}{2}\|U^{n+1}_h\|_{L^2}^2\biggl) \leq
C\epsilon^{-1}\|g^{n+1}\|_{L^2}^2
+\epsilon\|U^{n+1}_h\|_{L^2}^2,
\end{align*}
which further implies that (via Gronwall's inequality)
\begin{align}\label{d71}
&\max_{1\leq n\leq N}\|U^{n}_h\|_{L^2}\leq C .
\end{align}
Therefore,
\begin{align}\label{djt2}
\max_{1\leq n\leq N}\|U^{n}_h-u^{n}\|_{L^2}&\leq C\leq
\frac{C}{\max(\tau,h^{r+1})}(\tau+h^{r+1}) \leq
\frac{C}{\min(\tau_7,h^{r+1}_7)}(\tau+h^{r+1}).
\end{align}
Combining (\ref{proj-002}), (\ref{error-l2} ) and (\ref{djt2}), we see that
(\ref{optimalL2est}) holds unconditionally.

The proof of Theorem \ref{ErrestFEMSol} is completed. ~\endproof

\section{Proof of Theorem \ref{ErrestDisSol}}\label{sectD}
\setcounter{equation}{0}
First, we consider the Poisson equation 
\begin{equation} 
\left\{\begin{array}{ll}
-\Delta v=f -\frac{1}{|\Omega|}\int_\Omega f\d x  ,&\mbox{in}~~\Omega,\\
\partial_{\vec n}v=0 &\mbox{on}~~\partial\Omega,
\end{array}\right.
\label{poisson} 
\end{equation} 
in a convex polygon, and 
introduce some lemmas concerning the
$W^{2,p}$ and $H^{2+s}$ estimates of its solution. 
\begin{lemma}\label{lem00}
{\it Let $v$ be the solution of {\rm\refe{poisson}}
and $w \in W^{1,3}$ and $w_{\min}\leq w(x)\leq w_{\max}$,
where $w_{\min}$ and $w_{\max}$ are positive constants. 
If $f\in L^2$ and $\int_\Omega f\d x=0$, 
then $v \in H^2$ and for any $\epsilon\in(0,1/2)$ we have
\begin{align} 
&\|\nabla^2 v\|_{L^2}\leq \|f\|_{L^2}, \label{sj01}\\
&(1-\epsilon)\int_\Omega \sum_{i,j}|\partial_{ij} v |^2w\d x
\leq \int_\Omega |f|^2 w\d x
+C_{w_{\min},w_{\max},\|w\|_{W^{1,3}}}\epsilon^{-2}
\|\nabla v\|_{L^2}^{2}, \label{sj011}
\end{align}
}
\end{lemma}
\noindent{\it Proof}~~~ 
The inequality (\ref{sj01}) is a consequence of 
Theorem 3.1.1.1 in \cite{Grisvard}. 

For the inequality \refe{sj011}, we only present a priori estimates here.
By noting the identity
\begin{align*}
\partial_{ii}v\partial_{jj}v= \partial_i(\partial_iv\partial_{jj}v)
-\partial_j(\partial_iv\partial_{ij}v)+|\partial_{ij} v|^2, 
\end{align*}
we have 
\begin{align*}
\int_\Omega \sum_{i,j}|\partial_{ij}v|^2w\d x
&\leq \int_\Omega |f|^2 w\d x
+\int_{\Omega}\biggl(f\nabla w\cdot\nabla v 
-\sum_{i,j}\partial_jw\partial_i v\partial_{ij} v \biggl)\d x \,   
\end{align*}
and therefore,
\begin{align*} 
&(1-\epsilon)\int_\Omega \sum_{i,j}|\partial_{ij}v|^2w\d x\\
&\leq (1+\epsilon)\int_\Omega |f|^2 w\d x
+C\epsilon^{-1}\|w^{-1/2}\nabla w\|_{L^3}^2\|\nabla v\|_{L^6}^2\\
&\leq (1+\epsilon)\int_\Omega |f|^2 w\d x
+C\epsilon^{-1}(\|\nabla v\|_{L^2}^{2}+\|\nabla v\|_{L^2}^{4/3}\|\nabla^2 v\|_{L^2}^{2/3})\\
&\leq (1+\epsilon)\int_\Omega |f|^2 w\d x
+\epsilon \int_\Omega \sum_{i,j}|\partial_{ij} v|^2w\d x
+C\epsilon^{-2} \|\nabla v\|_{L^2}^{2}.
\end{align*}
\refe{sj011} follows immediately. ~\endproof

It can be found in literatures, such as Theorem 4.3.2.3 and Theorem 4.4.3.7 of 
\cite{Grisvard}, and (23.3) of \cite{Dauge},
that \begin{align}
&\|\nabla^2 v\|_{L^{1+p_*/2}}\leq C_*\|f\|_{L^{1+p_*/2}},\label{ddu}\\
&\|\nabla^2 v\|_{H^{s_*/2}}\leq C_*\|f\|_{L^{s_*/2}}
\label{ddu2}
\end{align}
for some positive constant $C_*\geq 4$, where
$p_*=\min(5/2,1/[1-\pi/(2\omega_m)])$, $s_*=\pi/\omega_m -1$ and
$\omega_m$ denotes the maximal interior angle
of the convex polygon $\Omega$. Since 
$\big\|f -\frac{1}{|\Omega|}\int_\Omega f\d x\big\|_{L^2}\leq \|f\|_{L^2}$,
the operator from 
$f$ to $\nabla^2 v$ defined by \refe{poisson} 
satisfies \refe{sj01} and \refe{ddu}-\refe{ddu2}.
By applying the complex interpolation 
(see Theorem 5.6.3 of \cite{BL}) 
to \refe{sj01} and \refe{ddu}-\refe{ddu2}, we obtain the following lemma. 

\begin{lemma}\label{lem001}
{\it 
Assume that $v \in H^2(\Omega)$ is the solution of the
equation {\rm\refe{poisson}}.  
Then 
\begin{align}\label{sj012}
&\|\nabla^2 v\|_{L^p}
\leq (1+\varepsilon_p) \|f\|_{L^p}\\
&\|\nabla^2 v\|_{H^s}
\leq (1+\overline\varepsilon_s) \|f\|_{H^s}
\label{lemH2s}
\end{align}
for $p\in(2,p_*)$ and $s\in(0,s_*)$, where 
$\lim_{p\rightarrow 2}\varepsilon_p=0$ and 
$\lim_{s\rightarrow 0}\overline\varepsilon_s=0$.
}
\end{lemma}

Based on the regularity assumption (\ref{regularity}), we set
\begin{align*}
& K=\|u\|_{L^\infty(\Omega \times (0,T))}+\|\nabla
u\|_{L^\infty(\Omega \times (0,T))}+2 .
\end{align*}
Then, by the regularity assumptions on $\sigma$, there exist
positive constants $0<\sigma_K<1$ and $C_K$ such that for
$0\leq s\leq K$ we
have
\begin{align}
&\sigma(s^2)\geq \sigma_K ,~~~
|\sigma(s^2)|+|\sigma'(s^2)|+|\sigma''(s^2)|\leq C_K ,
\label{g-sigma}
\end{align}
and we choose $p$ so close to $2$ that 
\begin{align}
\varepsilon_p< \lambda^2\sigma_K^{2}.
\label{pp}
\end{align}

Now we start to prove Theorem \ref{ErrestDisSol}. 
For the given $U^n\in H^{2+s_n}$, \refe{e-TDparab-1} can be viewed as a linear
elliptic boundary value problem and therefore, 
it 
admits a unique solution $U^{n+1}\in H^{2+s_{n+1}}$ for some positive constant 
$s_{n+1}>0$ (a qualitative regularity as a consequence of Lemma 4.2).   
Here we only prove the quantitative estimates
\refe{dt-bound}-\refe{dt-err}.

Before we study the error estimates \refe{optimalL2est},  
we prove by mathematical induction the following inequalities 
\begin{align}
&\|U^n \|_{L^\infty} +\| \nabla U^n \|_{L^\infty} \leq K ,
\label{a-1}\\
&\| e^n
\|_{W^{2,p}} \leq \tau^{1/3} \label{a-11}
\end{align}
assuming $\tau<\tau_0^*$ for some $\tau^*_0 >0$. 
Since $U^0=u_0$, the above inequalities hold for $n=0$.
We assume that (\ref{a-1})-(\ref{a-11}) hold for $0\leq n\leq k$ for
some nonnegative integer $k$, and
 prove the inequalities for $n=k+1$.

From \refe{e-parab-1}-\refe{IniC} and
\refe{e-TDparab-1}, we see that
$e^{n+1}$ satisfies the equation
\begin{align}
&D_\tau e^{n+1}-\nabla\cdot(\sigma(|\nabla u^n|^2)\nabla
e^{n+1})
\label{err-TDparab-01}\\
&={\cal E}^{n+1}_{\rm tr}-\nabla\cdot((\sigma(|\nabla
U^n|^2)-\sigma(|\nabla u^n|^2))\nabla U^{n+1}) , \nn
\end{align}
with the boundary condition $\nabla e^{n+1}\cdot{\vec n}=0$ 
and the initial condition $e^0=0$, where
\begin{align*}
{\cal E}^{n+1}_{\rm tr}=\partial_tu^{n+1} -D_\tau
u^{n+1}+\nabla\cdot [(\sigma(|\nabla u^n|^2)-\sigma(|\nabla u^{n+1}|^2))\nabla u^{n+1}]
\end{align*}
is the truncation error due to the time
discretization. By the regularity assumption
(\ref{regularity}), we have
\begin{align}
\max_{1\leq n\leq N}\|{\cal E}^{n}_{\rm tr}\|_{L^2}\leq
C,\qquad
\sum_{n=1}^{N}\tau \|{\cal E}^{n}_{\rm tr}\|_{L^2}^2\leq
C\tau^2 .
\end{align}
With a similar approach to \refe{sigma-h-1},  we can derive that
\begin{align}
&(\sigma(|\nabla U^n|^2) -\sigma(|\nabla u^n|^2)) \nabla
U^{n+1}\nn\\
&= \big(-2 \sigma'(|\nabla u^n|^2) \nabla u^n \cdot \nabla
e^n +
\sigma'(|\nabla u^n|^2)|\nabla e^n|^2 \big)\nabla U^{n+1}\nn
\\
&~~~ +\frac{1}{2}\sigma''(\xi^n) |(\nabla u^n+\nabla
U^n)\nabla
e^n|^2 \nabla U^{n+1}
\nn \\
& =-2 \sigma'(|\nabla u^n|^2) (\nabla u^n \cdot\nabla e^n)
\nabla
u^n
\nn \\
&~~~ - 2 \sigma'(|\nabla u^n|^2) (\nabla u^n \cdot\nabla
e^n) (\tau
D_{\tau} \nabla u^{n+1} - \nabla e^{n+1}) \nn\\
&~~~+\Big(\sigma'(|\nabla u^n|^2)|\nabla
e^n|^2+\frac{1}{2}\sigma''(\xi^n) |(\nabla u^n+\nabla
U^n)\nabla
e^n|^2\Big)\cdot\nabla U^{n+1}\nn\\
&\leq \gamma(|\nabla u^n|^2)|\nabla e^n|+ C \tau |\nabla
e^n| +
C|\nabla e^n||\nabla e^{n+1}| +C |\nabla e^n|^2
,\label{sigma-2}
\end{align}
where $\gamma(\cdot)$ is defined in \refe{gamma}.

Multiplying (\ref{err-TDparab-01}) by $e^{n+1}$ and using
(\ref{sigma-2}), we get
\begin{align*}
& D_\tau\biggl(\frac{1}{2}\|e^{n+1}\|_{L^2}^2\biggl) + \|
\sqrt{\sigma(|\nabla u^n|^2)} \nabla e^{n+1} \|_{L^2}^2 \\
&\leq \frac{1}{2} \| \sqrt{\gamma(|\nabla u^n|^2)} \nabla
e^n\|_{L^2}^2+\frac{1}{2}\| \sqrt{\gamma(|\nabla u^n|^2)}
\nabla
e^{n+1} \|_{L^2}^2 + C \tau(\| \nabla e^{n+1} \|_{L^2}^2 +
\|
\nabla e^{n}
\|_{L^2}^2)   \nn \\
&~~~+C \|\nabla e^n \|_{L^\infty}( \|\nabla e^{n}\|_{L^2}^2
+
\|\nabla e^{n+1}\|_{L^2}^2  )
+ \|{\cal E}^{n+1}_{\rm tr}\|_{L^2}^2 +\|e^{n+1}\|_{L^2}^2 ,
\end{align*}
which implies that
\begin{align}
& D_\tau\biggl(\frac{1}{2}\|e^{n+1}\|_{L^2}^2\biggl)+
\frac{1}{2}
\left ( \| \sqrt{\sigma(|\nabla u^n|^2)} \nabla e^{n+1}
\|_{L^2}^2
- \| \sqrt{\gamma(|\nabla u^n|^2)} \nabla e^{n} \|_{L^2}^2
\right ) \nn \\
&\leq C\tau^{1/4}( \|\nabla e^{n+1} \|_{L^2}^2 + \| \nabla
e^{n}
\|_{L^2}^2) + C \|e^{n+1}\|_{L^2}^2 +C\|{\cal E}_{\rm
tr}^{n+1}\|_{L^2}^2  ,\label{err-2}
\end{align}
where we have used (\ref{a-11}). By noting 
\begin{align*}
& \| \sqrt{\sigma(|\nabla u^n|^2)} \nabla e^{n+1} \|_{L^2}^2
- \|
\sqrt{\gamma(|\nabla u^n|^2)} \nabla e^{n} \|_{L^2}^2
\\
&=\| \sqrt{\sigma(|\nabla u^n|^2)-\gamma(|\nabla u^n|^2)}
\nabla
e^{n+1} \|_{L^2}^2 +\| \sqrt{\gamma(|\nabla u^{n}|^2)}
\nabla
e^{n+1} \|_{L^2}^2- \| \sqrt{\gamma(|\nabla u^{n-1}|^2)}
\nabla
e^{n} \|_{L^2}^2
\nn\\
&~~~ - \big((\gamma(|\nabla u^n|^2)-\gamma(|\nabla
u^{n-1}|^2)) \nabla e^{n},\nabla e^{n}\big) \\
& \geq \|\lambda\sigma(|\nabla u^n|^2)^{3/2}
\nabla
e^{n+1} \|_{L^2}^2+\tau D_\tau\| \sqrt{\gamma(|\nabla
u^n|^2)}
\nabla e^{n+1} \|_{L^2}^2 - \tau \| \sqrt{|D_{\tau}
\gamma(|\nabla
u^n|^2)|} \nabla e^{n} \|_{L^2}^2 \\
& \geq \lambda^{2}\sigma_K^3\|\nabla
e^{n+1} \|_{L^2}^2+\tau D_\tau\| \sqrt{\gamma(|\nabla
u^n|^2)}
\nabla e^{n+1} \|_{L^2}^2 - C\tau \|\nabla e^{n} \|_{L^2}^2 ,
\end{align*}
(\ref{err-2}) reduces to
\begin{align*}
& D_\tau\biggl(\frac{1}{2}\|e^{n+1}\|_{L^2}^2 +\frac{\tau}{2}\|
\sqrt{\gamma(|\nabla u^n|^2)} \nabla e^{n+1} \|_{L^2}^2\biggl)+
\frac{\lambda^2\sigma_K^3}{2} \| \nabla e^{n+1} \|_{L^2}^2 \nn \\
&\leq C\tau^{1/4}(\|\nabla e^{n+1} \|_{L^2}^2 + \|\nabla e^{n} \|_{L^2}^2 )
+ C\|e^{n+1}\|_{L^2}^2 +C\|{\cal E}_{\rm tr}^{n+1}\|_{L^2}^2  .
\end{align*}
By Gronwall's inequality, when $\tau<\tau_8$ for some positive constant $\tau_8$, we have
\begin{align}
\label{dt-l2} &
\max_{0\leq n\leq k}\|
e^{n+1}\|_{L^2}^2+\sum_{n=0}^k\tau\|e^{n+1}\|_{H^1}^2\leq
C\sum_{n=0}^{k}\tau \|{\cal E}_{\rm tr}^{n+1}\|_{L^2}^2\leq
C\tau^2 .
\end{align}
From the above inequality we also see that
\begin{align}
\label{dt-b}
&\|U^{n+1}\|_{L^2}\leq \| u^{n+1}\|_{L^2}+\| e^{n+1}\|_{L^2}
\leq C,\\
&\|D_\tau U^{n+1}\|_{L^2}\leq \|D_\tau
u^{n+1}\|_{L^2}+\|D_\tau
e^{n+1}\|_{L^2}  \leq C .
\end{align}
We rewrite (\ref{err-TDparab-01}) as
\begin{align}\label{dk29}
&D_\tau e^{n+1}-\sigma(|\nabla u^n|^2)\Delta
e^{n+1} \nn\\
&={\cal E}^{n+1}_{\rm tr}+2\sigma'(|\nabla U^n|^2)(\nabla^2U^n\nabla
U^{n})\cdot
\nabla e^{n+1} -(\sigma(|\nabla
U^n|^2)-\sigma(|\nabla u^n|^2)\Delta u^{n+1}) \nn\\
&~~~-\big[2\sigma'(|\nabla
U^n|^2)\nabla^2U^n\nabla
U^{n}-2\sigma'(|\nabla u^n|^2)\nabla^2u^n\nabla
u^{n}\big]\cdot\nabla u^{n+1} \nn\\
&~~~-[\sigma(|\nabla u^n|^2)-\sigma(|\nabla U^n|^2)]\Delta
e^{n+1} \nn\\
&={\cal E}^{n+1}_{\rm tr}+2\sigma'(|\nabla U^n|^2)(\nabla^2U^n\nabla
U^{n})\cdot
\nabla e^{n+1} -(\sigma(|\nabla
U^n|^2)-\sigma(|\nabla u^n|^2)\Delta u^{n+1}) \nn\\
&~~~+2\sigma'(|\nabla
u^n|^2)\nabla^2e^n\nabla
u^{n}\cdot\nabla u^{n}-[2\sigma'(|\nabla
U^n|^2)\nabla^2u^n\nabla
U^{n}-2\sigma'(|\nabla u^n|^2)\nabla^2u^n\nabla
u^{n}]\cdot\nabla u^{n+1} \nn\\
&~~~+[2\sigma'(|\nabla
U^n|^2)\nabla
U^{n}-2\sigma'(|\nabla
u^n|^2)\nabla
u^{n}]\cdot(\nabla^2e^n\nabla u^{n+1}) \nn\\
&~~~+2\tau \sigma'(|\nabla
u^n|^2)\nabla^2e^n\nabla
u^{n}\cdot\nabla D_\tau u^{n+1}
-[\sigma(|\nabla u^n|^2)-\sigma(|\nabla U^n|^2)]\Delta
e^{n+1}  .
\end{align}
Multiplying the above equation by $-\Delta e^{n+1}$ leads to 
\begin{align*}
&D_\tau \biggl(\frac{1}{2}|\nabla e^{n+1}|^2\d x\biggl)
+ \int_\Omega\sigma(|\nabla u^n|^2)|\Delta e^{n+1}|^2\d x\\
&\leq \|{\cal E}^{n+1}_{\rm tr}\|_{L^2}\|\Delta e^{n+1}\|_{L^2}
+ C\|\nabla^2U^n\|_{L^p}\|\nabla
U^{n}\|_{L^\infty}\|
\nabla e^{n+1}\|_{L^{2p/(p-2)}} \|\Delta e^{n+1}\|_{L^2}\\
&+C\|
\nabla e^{n+1}\|_{L^2} \|\Delta e^{n+1}\|_{L^2}
+\int_\Omega \gamma(|\nabla u^n|^2)
|\nabla^2e^{n}||\Delta e^{n+1}|\d x
+C\|\nabla e^n\|_{L^2}\|\Delta e^{n+1}\|_{L^2}\\
&~~~+C\|\nabla e^n\|_{L^\infty}\|\nabla^2 e^n\|_{L^2}
\|\Delta e^{n+1}\|_{L^2}
+C\tau\|\nabla^2 e^n\|_{L^2}\|\Delta e^{n+1}\|_{L^2}
+C\|\nabla e^n\|_{L^\infty}\|\Delta e^{n+1}\|_{L^2}^2\\
&\leq \frac{1}{2} \int_\Omega (\gamma(|\nabla u^n|^2)
+\epsilon+C\tau^{1/4})|\nabla^2e^{n}|^2\d x
+\frac{1}{2} \int_\Omega\gamma(|\nabla u^n|^2)|\Delta e^{n+1}|^2\d x\\
&~~~+C\epsilon^{-1}(\|{\cal E}^{n+1}_{\rm tr}\|_{L^2}^2
+\|\nabla e^n\|_{L^2}^2) ,
\end{align*}
which further reduces to
\begin{align*}
&D_\tau \biggl(\frac{1}{2}|\nabla e^{n+1}|^2\d x\biggl)
+ \frac{1}{2}\int_\Omega\sigma(|\nabla u^{n+1}|^2)
|\Delta e^{n+1}|^2\d x
+ \frac{\lambda^2}{2}\int_\Omega\sigma(|\nabla u^n|^2)^3
|\Delta e^{n+1}|^2\d x\\
&\leq \frac{1}{2} \int_\Omega (\gamma(|\nabla u^n|^2)
+\epsilon+C\tau^{1/4})|\nabla^2e^{n}|^2\d x
+C\epsilon^{-1}(\|{\cal E}^{n+1}_{\rm tr}\|_{L^2}^2
+\|\nabla e^n\|_{L^2}^2)\\
&~~~+ \frac{1}{2}\int_\Omega\big[\sigma(|\nabla u^{n+1}|^2)
-\sigma(|\nabla u^{n}|^2)\big]
|\Delta e^{n+1}|^2\d x\\
&\leq \frac{1}{2} \int_\Omega (\gamma(|\nabla u^n|^2)
+\epsilon+C\tau^{1/4})|\nabla^2e^{n}|^2\d x
+ C\tau\int_\Omega|\Delta e^{n+1}|^2\d x
+C\epsilon^{-1}(\|{\cal E}^{n+1}_{\rm tr}\|_{L^2}^2
+\|\nabla e^n\|_{L^2}^2) .
\end{align*}
When $\tau<\tau_9$ for some positive constant $\tau_9$, we get
\begin{align*}
&D_\tau \biggl(\frac{1}{2}|\nabla e^{n+1}|^2\d x\biggl)
+ \frac{1}{2}\int_\Omega\Big[\sigma(|\nabla u^{n+1}|^2)
+\frac{\lambda^2}{2}\sigma(|\nabla u^n|^2)^3\Big]
|\Delta e^{n+1}|^2\d x\\
&\leq \frac{1}{2} \int_\Omega (\gamma(|\nabla u^n|^2)
+2\epsilon+C\tau^{1/4})|\nabla^2e^{n}|^2\d x
+C\epsilon^{-1}(\|{\cal E}^{n+1}_{\rm tr}\|_{L^2}^2
+\|\nabla e^n\|_{L^2}^2) .
\end{align*}
and by applying Lemma \ref{lem00} with
$w=\sigma(|\nabla u^{n+1}|^2)
+\frac{\lambda^2}{2}\sigma(|\nabla u^{n}|^2)^3$,
we obtain
\begin{align*}
&D_\tau \biggl(\frac{1}{2}|\nabla e^{n+1}|^2\d x\biggl)
+ \frac{1}{2}\int_\Omega\Big[(1-\epsilon)\sigma(|\nabla u^{n+1}|^2)
+\frac{(1-\epsilon)\lambda^2}{2}\sigma(|\nabla u^n|^2)^3\Big]
|\nabla^2 e^{n+1}|^2\d x\\
&\leq \frac{1}{2} \int_\Omega [ \gamma(|\nabla u^n|^2)
+\epsilon
+C\tau^{1/4}]|\nabla^2 e^{n}|^2\d x
+C\epsilon^{-2}(\|{\cal E}^{n+1}_{\rm tr}\|_{L^2}^2
+\|\nabla e^n\|_{L^2}^2+\|\nabla e^{n+1}\|_{L^2}^2 )  .
\end{align*}
By choosing $\epsilon$ small enough and when
$\tau <\tau_{10}$ for some positive constant $\tau_{10}$, we derive that
\begin{align*}
&D_\tau \biggl(\frac{1}{2}|\nabla e^{n+1}|^2\d x\biggl)
+ \frac{1}{2}\int_\Omega\sigma(|\nabla u^{n+1}|^2)
|\nabla^2 e^{n+1}|^2\d x\\
&\leq \frac{1}{2}\int_\Omega\big(\sigma(|\nabla u^n|^2)
-\lambda^2\sigma_K^3/2\big)
|\nabla^2 e^{n}|^2\d x+C
\big(\|{\cal E}^{n+1}_{\rm tr}\|_{L^2}^2
+\|\nabla e^n\|_{L^2}^2+\|\nabla e^{n+1}\|_{L^2}^2\big) ,
\end{align*}
which in turn shows that (with Gronwall's inequality)
\begin{align}\label{ly6}
&\max_{0\leq n\leq k}\frac{1}{2}\|\nabla e^{n+1}\|^2_{L^2}
+ \frac{\lambda^2\sigma_K^3}{8}
\sum_{n=1}^k\tau \|\nabla^2 e^{n+1}\|^2_{L^2}
\leq  C\sum_{n=1}^k\tau 
\|{\cal E}^{n+1}_{\rm tr}\|_{L^2}^2 \leq C\tau^2 .
\end{align}
From the above inequality we further derive that
\begin{align}\label{ly7}
&\max_{0\leq n\leq k} (\|U^{n+1}\|^2_{H^1}
+\|D_\tau U^{n+1}\|^2_{H^1})
+ \sum_{n=1}^k\tau \|D_\tau U^{n+1}\|^2_{H^2}
\leq  C .
\end{align}

From (\ref{dk29}) we see that
\begin{align*}
\|D_\tau e^{n+1}\|_{L^2}\leq C\|{\cal E}^{n+1}_{\rm tr}\|_{L^2}+C\|
e^{n}\|_{H^2}+C\|
e^{n+1}\|_{H^2},
\end{align*}
and by using (\ref{ly6}), 
\begin{align}
\sum_{n=0}^{k}\tau\|D_\tau
e^{n+1}\|_{L^2}^2\leq C\sum_{n=0}^{k}
\tau\|{\cal E}^{n+1}_{\rm tr}\|_{L^2}+C\sum_{n=0}^{k}\tau\|
e^{n+1}\|_{H^2}\leq C\tau^2.
\end{align}
In particular, the above inequality implies that
$\|D_\tau e^{k+1}\|_{L^2}\leq C\tau^{1/2}$
and $\|D_\tau e^{k+1}\|_{H^1}\leq C$  
from (\ref{ly7}). 
By an interpolation between $L^2$ and $H^1$, we have 
\begin{align*}
&\|D_\tau e^{k+1}\|_{L^p}\leq
C\|D_\tau e^{k+1}\|_{L^2}^{2/p}
\|D_\tau e^{k+1}\|_{H^1}^{1-2/p}\leq C\tau^{1/p}  .
\end{align*}
We rewrite (\ref{dk29}) by 
\begin{align} \label{dk30}
&\Delta
e^{k+1}=\sigma(|\nabla u^k|^2)^{-1}2\sigma'(|\nabla
u^k|^2)\nabla^2e^k\nabla
u^{k}\cdot\nabla u^{k}+G ,
\end{align}
where
\begin{align}\label{lkd7}
\|G\|_{L^p}&\leq C\|D_\tau e^{k+1}\|_{L^p}
+C\|{\cal E}^{k+1}_{\rm tr}\|_{L^p}
+(C\|\nabla^2U^k\|_{L^p}+C)\|\nabla e^{k+1}\|_{L^\infty} \nn\\
&~~~+C\tau \|\nabla^2e^k\|_{L^p}+C (\|\nabla^2e^{k}\|_{L^p}
+\|\nabla^2e^{k+1}\|_{L^p})\|\nabla e^{k}\|_{L^\infty} \nn\\
&\leq C\tau^{1/p}+\epsilon\|\nabla^2 e^{k+1}\|_{L^p}
+C_{\epsilon^{-1}}\|\nabla e^{k+1}\|_{L^2}
+C\tau^{1/4}(\|\nabla^2 e^{k+1}\|_{L^p}
+\|\nabla^2 e^{k}\|_{L^p}) \nn\\
&\leq C_{\epsilon^{-1}}\tau^{1/p}
+\epsilon(\|\nabla^2 e^{k+1}\|_{L^p}
+\|\nabla^2 e^{k}\|_{L^p}) .
\end{align}
With \refe{pp}, we apply \refe{sj012} to the elliptic equation
(\ref{dk30}) to get
\begin{align*}
\|\nabla^2 e^{k+1}\|_{L^p}&\leq (1+\lambda^2\sigma_K^{2})
\|\sigma(|\nabla
u^k|^2)^{-1}\gamma(|\nabla
u^k|^2)\nabla^2e^k\|_{L^p}+ (1+\lambda^2\sigma_K^{2})\|G\|_{L^p} \\
&\leq (1-\lambda^4\sigma_K^{4})\|\nabla^2e^k\|_{L^p}
+ (1+\lambda^2\sigma_K^{2})\|G\|_{L^p} .
\end{align*}
With 
$\epsilon=\lambda^4\sigma_K^4/(4+2\lambda^2\sigma_K^2)$ in (\ref{lkd7}), 
a straightforward calculation gives
\begin{align*}
\|\nabla^2 e^{k+1}\|_{L^p}
&\leq (1-\lambda^4\sigma_K^{4}/2)\|\nabla^2 e^{k}\|_{L^p}
+C\tau^{1/p}
\end{align*}
when  $\tau<\tau_{11}$ for some positive constant $\tau_{11}$.
By the Sobolev embedding inequality, we have
$\|e^k\|_{L^p} +\|\nabla e^k\|_{L^p}\leq C\|e^k\|_{H^2}\leq C\tau^{1/2}$
and therefore,
\begin{align}\label{lkq6}
\|e^{k+1}\|_{W^{2,p}}
&\leq (1-\lambda^4\sigma_K^{4}/2)\|e^{k}\|_{W^{2,p}}+C\tau^{1/p} 
\end{align}
which, by noting $1 < p \leq 5/2$, leads to 
\begin{align}\label{ew2pest}
\|e^{k+1}\|_{W^{2,p}}
&\leq \tau^{1/3} 
\end{align}
when $\tau<\tau_{12}$ for some positive constant $\tau_{12}$.
By using the Sobolev embedding inequality again, we obtain
\begin{align*}
\|e^{k+1}\|_{L^\infty} +\| \nabla e^{k+1}\|_{L^\infty}
\leq C\|e^{k+1}\|_{W^{2,p}}\leq C\tau^{1/3}
\end{align*}
which further implies that
\begin{align} \label{d27h}
&\|e^{k+1}\|_{L^\infty} +\| \nabla e^{k+1}\|_{L^\infty}
 \leq \tau^{1/4} ,\\
&\|U^{k+1}\|_{L^\infty} +\| \nabla U^{k+1}\|_{L^\infty} \leq K ,
\end{align}
when $\tau< \tau_{13}$ for some positive constant $\tau_{13}$.

The induction on (\ref{a-1})-(\ref{a-11}) is closed,
and \refe{dt-l2} and \refe{ly6}-\refe{ew2pest}
hold for $k=N$ provided $\displaystyle\tau<\tau_{0}^*:
=\min_{8\leq i\leq 13}\tau_i$.

It remains to estimate $\|U^{n+1}\|_{H^{2+s}}$ for some $s>0$. From \refe{ew2pest} we see that $\nabla U\in C^\alpha $ for some $\alpha>0$. Rewrite \refe{e-TDparab-1} as 
\begin{align}\label{DeltaUn}
 -\Delta U^{n+1} 
&=\frac{\sigma'(|\nabla U^n|^2)}{\sigma(|\nabla U^n|^2)} 
(\nabla^2 U^n\nabla U^n)\cdot\nabla U^{n} \nn\\
&~~~+\tau\frac{\sigma'(|\nabla U^n|^2)}{\sigma(|\nabla U^n|^2)} 
(\nabla^2 U^n\nabla U^n)\cdot\nabla D_\tau U^{n+1} +g^{n+1}-D_\tau U^{n+1} \nn\\
&=l(\nabla^2 U^n)+g^{n+1}-D_\tau U^{n+1} ,
\end{align}
where the linear operator $l$ defined by
$$ 
l(\nabla^2U^n)= \frac{\sigma'(|\nabla U^n|^2)}{\sigma(|\nabla U^n|^2)} 
(\nabla^2 U^n\nabla U^n)\cdot\nabla U^{n}  
+\tau\frac{\sigma'(|\nabla U^n|^2)}{\sigma(|\nabla U^n|^2)} 
(\nabla^2 U^n\nabla U^n)\cdot\nabla D_\tau U^{n+1}  
$$
satisfies that
\begin{align}
&\|l(\nabla^2U^n)\|_{L^2}\leq \bigg(\frac{K^2}{\lambda^2+K^2}+C\tau^{1/3}\bigg) 
\|\nabla^2U^n\|_{L^2} \nn\\[5pt]
&\|l(\nabla^2U^n)\|_{H^\alpha}\leq C\|\nabla^2U^n\|_{H^\alpha} \nn .
\end{align}
By choosing $\tau$ small enough and using the complex interpolation 
between $L^2$ and $H^\alpha$ we derive that, there exist positive constants $s_K$ such that 
\begin{align}
\|l(\nabla^2U^n)\|_{H^{s}}
&\leq \bigg(\frac{K^2}{\lambda^2+K^2}
+C\tau^{1/3}\bigg)^{1-s_K/\alpha}C^{s_K/\alpha}\|\nabla^2U^n\|_{H^{s}}\nn\\
&\leq \bigg(1-\frac{\lambda^2}{2\lambda^2+2K^2}\bigg)\|\nabla^2U^n\|_{H^{s}} 
\quad\mbox{for}~~s\in[0,s_K] .
\end{align}
Therefore, by applying \refe{lemH2s} to the equation \refe{DeltaUn} 
we obtain that
\begin{align}
\|\nabla^2U^{n+1} \|_{H^{s}}\leq (1+\overline\varepsilon_s)\bigg[
\bigg(1-\frac{\lambda^2}{2\lambda^2+2K^2}\bigg)\|\nabla^2U^n\|_{H^{s}}
+\|g^{n+1}\|_{H^s}
+\|D_\tau U^{n+1}\|_{H^s}\bigg] ,
\end{align}
and choosing $s_0$ so small that $\overline\varepsilon_{s_0}<\lambda^2/(2\lambda^2+2K^2)$,  
we get
\begin{align}
\|\nabla^2U^{n+1} \|_{H^{s_0}}\leq
\bigg(1-\frac{\lambda^4}{4(\lambda^2+K^2)^2}\bigg)\|\nabla^2U^n\|_{H^{s_0}}
+C\|g^{n+1}\|_{H^{s_0}}
+C\|D_\tau U^{n+1}\|_{H^{s_0}} .
\end{align}
Iterations of the above inequality give 
\begin{align}
\max_{1\leq n\leq N}
\|\nabla^2U^{n} \|_{H^{s_0}}\leq  
C(\max_{1\leq n\leq N}\|g^{n}\|_{H^{s_0}}
+\max_{1\leq n\leq N}\|D_\tau U^{n}\|_{H^{s_0}} ) \leq C  .
\end{align}
The proof of Theorem \ref{ErrestDisSol} is completed.
~\endproof

\section{Numerical example}
\setcounter{equation}{0}
In this section, we present an example to confirm our theoretical analysis. 
All computations are performed by FreeFEM++ in double precision \cite{Hecht}.

We solve (\ref{e-parab-1})-(\ref{IniC}) in the domain
$\Omega= [0,1] \times [0, 1]$ up to the time $T=1$, where the diffusion coefficient 
$\sigma(|\nabla u|^2)$ is given by \refe{sigma}, the function
$g$ and
$u_0$ are chosen corresponding to the exact solution
\begin{align}\label{sd8}
&u(x,y, t) = e^{0.01t }\cos (2\pi x)\cos (2\pi y)/4 .
\end{align}
To test the convergence rate in the spatial direction, 
a uniform triangulation is generated with
$M+1$ points on each side of the rectangular domain with $h=\sqrt{2}/M$, 
and we choose a very small time step $\tau =2^{-15}$. In this case, 
the optimal error estimate given in Theorem \ref{ErrestFEMSol} 
is, approximately,
$$ 
\|U^n_h-u^n\|_{L^2} = O(h^{r+1}) \, . 
$$ 
For $\lambda=1$, 
we present the $L^2$-norm errors in Table \ref{Tab1}, where 
the convergence rate is calculated based on the
numerical results corresponding to two finer meshes.
We see that the $L^2$-norm errors are proportional to
$h^{r+1}$, which is consistent with our theoretical error analysis.

For comparison, we also present the numerical results for the case of $\lambda=0.2$ in Table \ref{Tab2}, with the quadratic FEM. We can 
see that the convergence of the numerical solution for the problem with 
$\lambda=0.2$ is much worse than the convergence of the numerical solution with $\lambda=1$. This indicates that our error estimate presented in this paper does not hold uniformly as $\lambda\rightarrow 0$.

\bigskip

\begin{table}[ht]
\vskip-0.2in
\begin{center}
\caption{$L^2$-norm errors of the numerical 
solution for $\lambda=1$}\vskip
0.1in
\label{Tab1}
\begin{tabular}{c|c|ccc}
\hline     $M$ & $\|U^N_h-u^N\|_{L^2}$~for~$r=2$
&$\|U^N_h-u^N\|_{L^2}$~for~$r=3$
      \\
\hline
 8 & 9.0361E-04 & 3.6292E-04   \\
 16 & 1.1846E-04 & 7.6558E-05   \\
 32 & 1.4948E-05 & 4.1758E-07    \\
\hline convergence rate & $O(h^{3.0})$
&  $O(h^{4.1})$     \\
\hline
\end{tabular}
\end{center}
\end{table}


\begin{table}[ht]
\vskip-0.2in
\begin{center}
\caption{$L^2$-norm errors of the numerical 
solution for $\lambda=0.2$}\vskip
0.1in
\label{Tab2}
\begin{tabular}{c|cccc}
\hline     $M$ & $\|U^N_h-u^N\|_{L^2}$~for~$r=2$

      \\
\hline
 8 &   5.3586E-02    \\
 16 & 1.0428E-02    \\
 32 & 2.7755E-04      \\
 64 & 9.0595E-06       \\
128 & 1.1281E-06     \\
\hline convergence rate &   $O(h^{3.0})$
       \\
\hline
\end{tabular}
\end{center}
\end{table}

\begin{figure}[h]
\begin{minipage}[b]{0.5\linewidth}
\centering
\includegraphics[width=\textwidth]{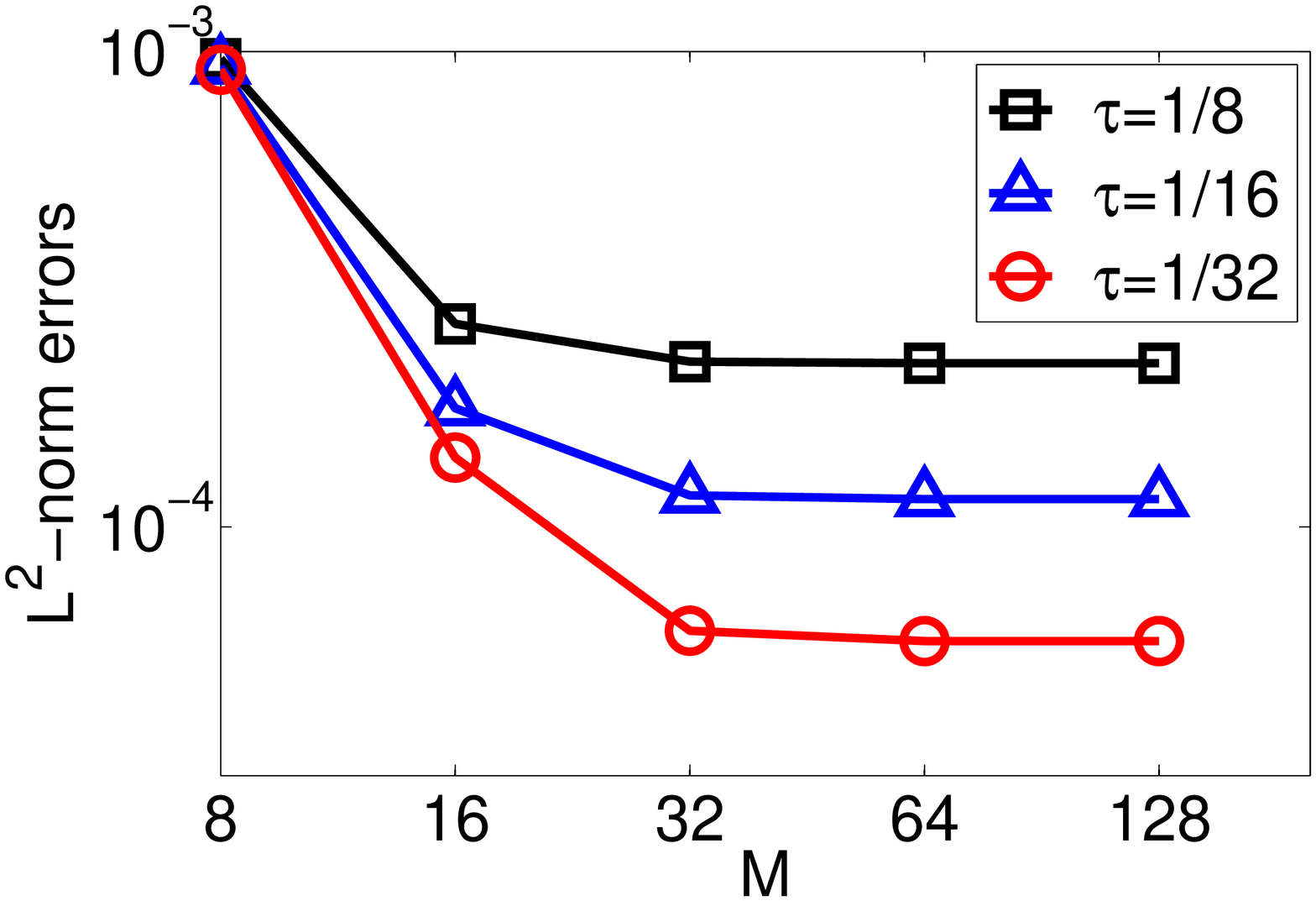}
\caption{$L^2$-norm error with $r=2$}\label{Fig1}
\end{minipage}
\hspace{-0.0cm}
\begin{minipage}[b]{0.5\linewidth}
\centering
\includegraphics[width=\textwidth]{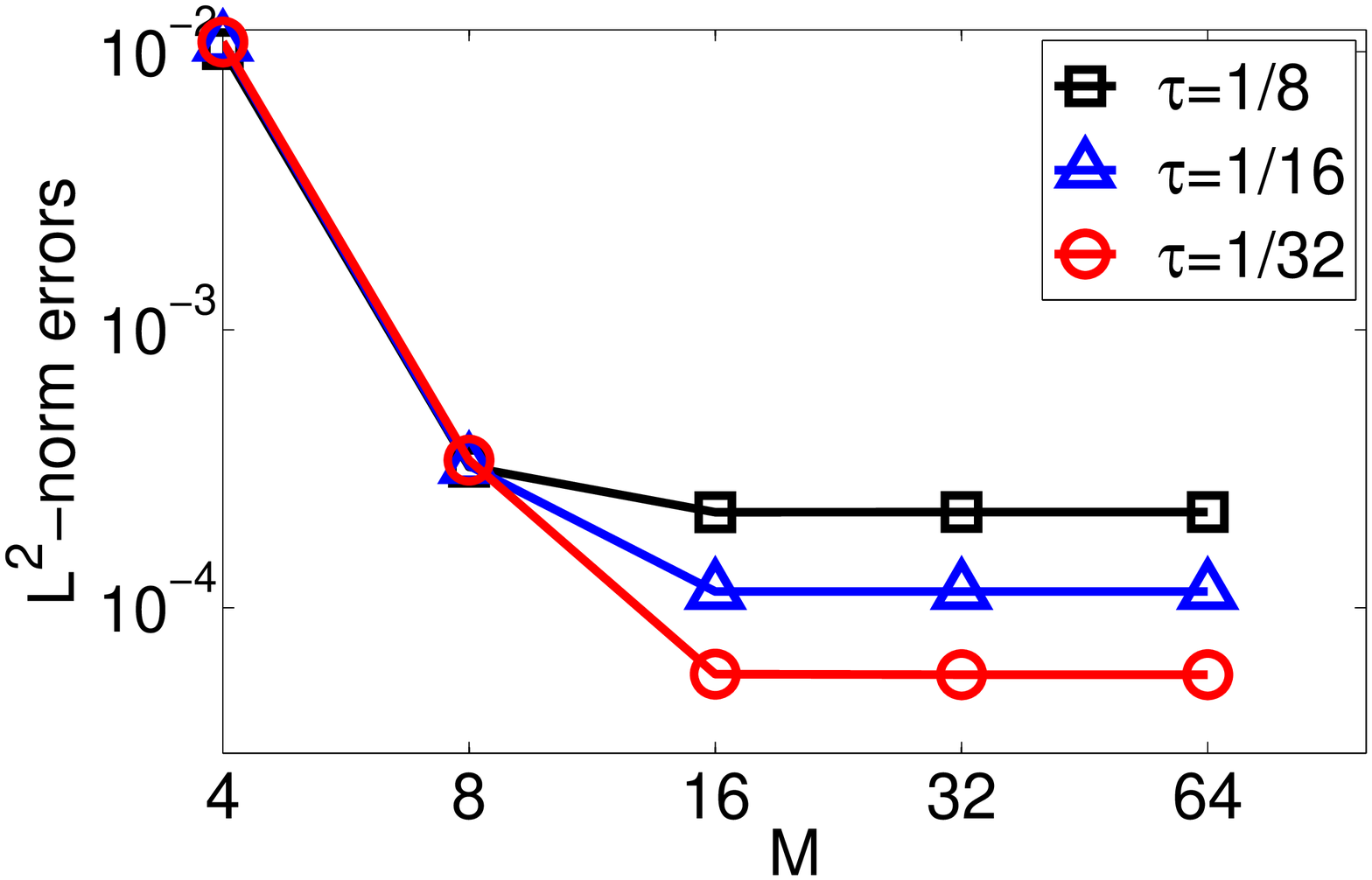}
\caption{$L^2$-norm error with $r=3$.}\label{Fig2}
\end{minipage}
\end{figure}

To test the convergence rate in the temporal direction and 
the stability of the numerical solution, we solve
(\ref{e-parab-1})-(\ref{IniC}) with several refined meshes 
for each fixed $\tau$. The $L^2$-norm errors of the
numerical solution are presented in Figure \ref{Fig1} and
Figure \ref{Fig2} for $r=2,3$, respectively, 
on the logarithmic scale. We see that, for each fixed $\tau$, 
the $L^2$-norm error of the numerical solution tends to a constant 
which is proportional to $\tau$. Therefore, no
restriction on the grid ratio is needed.

\section{Conclusion}
In this paper, we have presented optimal error estimates for a linearized
backward Euler--Galerkin FEM for a nonlinear and non-degenerate diffusion equation in
a convex polygonal domain based on certain assumption
on the regularity of the exact solution. 
For this strongly nonlinear equation,
no previous works have been devoted to the error analysis for linearized semi-implicit
FEMs, and existing analyses for implicit schemes still require 
certain restrictions on the time stepsize.
Our analysis shows that the numerical solution of the linearized semi-implicit scheme
achieves optimal convergence rate without any time-step condition. 

There are some applications in which 
some degenerate diffusion equations 
($\lambda =0$) should be investigated, 
such as total variation model \cite{BC,FOP,FP} and parabolic $p$-Laplacian \cite{BLiu,DER,VO} 
without regularization. 
Numerical analysis for such degenerate equations is extremely difficult. 
Existing techniques in classical FEMs may not work well. 
Analysis for linearized schemes was less explored and 
many efforts focused  only on implicit schemes and suboptimal error estimates 
due to the degeneracy. 
The extension of our analysis to the nonlinear non-degenerate diffusion equation 
in three-dimensional space and to the the nonlinear degenerate 
equations is our future works.

\bigskip

\section*{Appendix: $\bf H^{2+\delta_0}$ 
regularity of the equation $\bf \nabla\cdot(A(\nabla U^n)\nabla\psi)=\varphi$}
\renewcommand{\thelemma}{A.\arabic{lemma}}
\renewcommand{\theproposition}{A.\arabic{lemma}}
\renewcommand{\theequation}{A.\arabic{equation}}
\setcounter{lemma}{0} \setcounter{equation}{0}

Under the assumption that $\|U^n\|_{H^{2+s_0}}\leq C$ (as given in Theorem \ref{ErrestDisSol}) 
we consider the equation
\begin{align}
\left\{
\begin{array}{ll}
\nabla\cdot(A(\nabla U^n)\nabla\psi)=\varphi
&\mbox{in}~~\Omega ,\\[5pt]
\nabla\psi\cdot{\vec n}=0 &\mbox{on}~~\partial\Omega,
\end{array}\right.
\end{align}
with the compatibility condition $\int_\Omega\varphi\d x=0$ and the normalization 
condition $\int_\Omega\psi\d x=0$. 

For simplicity, we only present a priori estimates here. 
The equation can be written as 
\begin{align}\label{DepsiEq}
\Delta \psi 
= \frac{\nabla^2\psi \nabla U^n\cdot \nabla U^n }{\lambda^2
+|\nabla U^n|^2}
-\sqrt{1+|\nabla U^n|^2/\lambda^2}
\partial_iA_{ij}(\nabla U^n)\partial_{j}\psi 
+\sqrt{1+|\nabla U^n|^2/\lambda^2}\varphi  
\end{align}
and by using \refe{sj01} we derive that
\begin{align*}
\|\nabla^2\psi \|_{L^2}
&\leq 
\bigg\|\frac{| \nabla U^n|^2 }{\lambda^2+|\nabla U^n|^2}\bigg\|_{L^\infty}
\|\nabla^2\psi\|_{L^2}+C\|\nabla^2 U^n\|_{L^p}\|\nabla\psi \|_{L^{\bar p}}
+ C\|\varphi \|_{L^2}\\
&\leq 
\frac{K^2 }{1+K^2}
\|\nabla^2\psi\|_{L^2}+C\|\nabla\psi \|_{L^{\bar p}}
+ C\|\varphi \|_{L^2}, 
\end{align*}
which further reduces to
\begin{align*}
\|\nabla^2\psi \|_{L^2}
&\leq C\|\nabla\psi \|_{L^{\bar p}}
+ C\|\varphi \|_{L^2}\\
&\leq C\|\nabla\psi \|_{L^2}^{2/\bar p}\|\psi \|_{H^2}^{1-2/\bar p}
+ C\|\varphi \|_{L^2} \\
&\leq \epsilon\|\psi \|_{H^2} +C_\epsilon \|\nabla\psi \|_{L^2} + C\|\varphi \|_{L^2}\\
&\leq \epsilon\|\nabla^2\psi \|_{L^2} +C_\epsilon \|\varphi \|_{L^2} .
\end{align*}
From the last inequality we see that
\begin{align}
\|\psi \|_{H^2}
\leq C\|\varphi \|_{L^2} .
\end{align}

Since $\nabla U^n\in H^{1+s_0}\hookrightarrow C^{s_0}(\overline\Omega)$, it follows that
\begin{align*}
\begin{array}{ll}
\displaystyle\Big\|\frac{F\nabla U^n\cdot \nabla U^n }{\lambda^2
+|\nabla U^n|^2}\Big\|_{L^2}\leq 
\frac{K^2}{\lambda^2+K^2}\|F\|_{(L^2)^{2\times 2}} 
&\mbox{for}~~s\in(0,s_0),\\[15pt]
\displaystyle\Big\|\frac{F
\nabla U^n\cdot \nabla U^n }{\lambda^2+|\nabla U^n|^2}
\Big\|_{H^{s_0}}\leq C\|F\|_{(H^{s_0})^{2\times 2}} & \\[15pt]
\|\sqrt{1+|\nabla U^n|^2/\lambda^2}
\partial_iA_{ij}(\nabla U^n)\partial_{j}\psi \|_{H^s}\leq 
C\|\psi \|_{H^{1+s}} &\mbox{for}~~s\in(0,s_0) ,
\end{array}
\end{align*}
and by the complex interpolation method \cite{BL} we derive that there exists 
a positive constant $\overline s_0\in (0,\min (s_*,s_0))$ such that  
(where $s_*$ is given in \refe{lemH2s})
\begin{align*}
\displaystyle\bigg\|\frac{F
\nabla U^n\cdot \nabla U^n }{\lambda^2
+|\nabla U^n|^2}\bigg\|_{H^s}
&\leq 
\bigg(\frac{K^2}{\lambda^2+K^2}\bigg)^{1-\bar s_0/s_0}
C^{\bar s_0/s_0}\|F\|_{(H^s)^{2\times 2}}\\
&\leq\bigg(1-\frac{1}{2\lambda^2+2K^2}\bigg)
\|F\|_{(H^s)^{2\times 2}} &\mbox{for}~~s\in(0,\overline s_0).
\end{align*}
By applying \refe{lemH2s} to \refe{DepsiEq} 
and using the last inequality, we derive that
$$
\|\nabla^2\psi \|_{H^{s}}\leq 
 (1+\overline\varepsilon_s)  \bigg(1-\frac{1}{2\lambda^2+2K^2}\bigg)\|\nabla^2\psi 
\|_{H^{s}}+C\|\psi \|_{H^{1+s}}+C\|\varphi\|_{H^s}\quad\mbox{for}~~s\in(0, \overline s_0) .
$$
There exists a positive constant $\delta_0$ such that 
$ \overline\varepsilon_s<1/(2\lambda^2+2K^2)$ when 
$s\in(0, \delta_0]$, and the last inequality reduces to 
\begin{align*}
\|\nabla^2\psi \|_{H^{s}}\leq 
\bigg(1-\frac{1}{4(\lambda^2+K^2)^2}\bigg)\|\nabla^2\psi \|_{H^{s}} 
+C\|\psi \|_{H^{1+s}}+C\|\varphi\|_{H^{s}}\quad\mbox{for}~~s\in(0,\delta_0] .
\end{align*}
which further implies that
\begin{align}
\|\psi \|_{H^{2+s}}\leq 
C\|\varphi\|_{H^{s}}\quad\mbox{for}~~s\in(0,\delta_0] .
\end{align}

%

%

\end{document}